\documentclass{dk_article}
\usepackage[top=3cm, bottom=3cm, left=3cm, right=3cm]{geometry}
\usepackage{graphicx}
\usepackage{amsmath}
\usepackage{todonotes}
\usepackage{amsthm}
\usepackage{mhequ}
\usepackage{color}

\usepackage{calrsfs}
\DeclareMathAlphabet{\pazocal}{OMS}{zplm}{m}{n}

\usepackage{hyperref}
\usepackage{amssymb}
\usepackage{caption}
\usepackage{subcaption}

\definecolor{darkred}{rgb}{.7,0,0}

\definecolor{darkgreen}{rgb}{0,0.7,0}

\definecolor{darkblue}{rgb}{0,0,0.7}

\title{Fluctuations in the heterogeneous multiscale methods\\ for
  fast-slow systems}
\author{David Kelly and Eric Vanden-Eijnden}
\institute{\;Courant Institute, New York University, NY, USA.\\
  \email{\{dtkelly,eve2\}@cims.nyu.edu}\\[10pt]
{\it Dedicated with admiration and friendship to Bjorn Engquist on the occasion of his 70th birthday.}}

\newcommand{\Xbar}{\bar{X}}

\newcommand{\ftilde}{\tilde{f}}

\newcommand{\CS}{\mathcal{S}}

\newcommand{\BE}{\mathbf{E}}
\newcommand{\BP}{\mathbf{P}}

\newcommand{\xbar}{\bar{x}}

\newcommand{\CH}{\pazocal{H}}
\newcommand{\CHhat}{\widehat{\pazocal{H}}}

\newcommand{\zhat}{\hat{z}}
\newcommand{\Zhat}{\widehat{Z}}

\newcommand{\CL}{\mathcal{L}}

\newcommand{\CV}{\mathcal{V}}
\newcommand{\Ytilde}{\widetilde{Y}}
\newcommand{\Shat}{\widehat{S}}
\newcommand{\CN}{\pazocal{N}}
\newcommand{\Xtilde}{\widetilde{X}}

\begin{document}

\maketitle


\begin{abstract}
  How heterogeneous multiscale methods (HMM) handle fluctuations
  acting on the slow variables in fast-slow systems is
  investigated. In particular, it is shown via analysis of central
  limit theorems (CLT) and large deviation principles (LDP) that the
  standard version of HMM artificially amplifies these fluctuations. A
  simple modification of HMM, termed parallel HMM, is introduced and
  is shown to remedy this problem, capturing fluctuations correctly
  both at the level of the CLT and the LDP. Similar type of arguments
  can also be used to justify that the $\tau$-leaping method used in
  the context of Gillespie's stochastic simulation algorithm for
  Markov jump processes also captures the right CLT and LDP for these
  processes.
\end{abstract}

\large

 \section{Introduction}
 \label{s:intro}
 The heterogeneous multiscale methods (HMM)
 \cite{engquist03,vanden01,%
   weinan2007heterogeneous,abdulle2012heterogeneous} provide an
 efficient strategy for integrating fast-slow systems of the type
\begin{equs}
    \label{e:fs_intro}
    \frac{dX^\eps}{dt} &= f(X^\eps,Y^\eps) \\
    \frac{dY^\eps}{dt} &= \frac{1}{\eps} g(X^\eps,Y^\eps) \;.
\end{equs}
The method relies on an averaging principle that holds under some
assumption of ergodicity and states that as $\eps\to0$ the slow
variables $X^\eps$ can be uniformly approximated by the solution to
the following averaged equation
 \begin{equation}
   \label{eq:2}
   \dot{\Xbar} = F(\Xbar)\;.
\end{equation}
Here $F(x) = \int f(x,y) \mu_x(dy)$ is the averaged vector field, with
$\mu_x(dy)$ being the ergodic invariant measure of the fast variables
$Y_x$ with a frozen $x$ variable. This averaging principle is akin to
the law of large number (LLN) in the present context and it suggests
to simulate the evolution of the slow variables using~\eqref{eq:2}
rather than~\eqref{e:fs_intro} when $\eps$ is small. This requires to
estimate $F(x)$, which typically has to be done on-the-fly given the
current value of the slow variables. To this end, note that if Euler's
method with time step $\Delta t$ is used as integrator for the slow
variables in~\eqref{e:fs_intro}, we can approximate
$X^\eps(n\Delta t)$ by $x_n$ satisfying the recurrence
\begin{equ}
  \label{e:scheme1}
  x^\eps_{n+1} = x^\eps_{n} + \int_{n\Delta t}^{(n+1)\Delta t}
  f(x_n^\eps,Y^\eps_{x_n^\eps}(s))ds \;,
\end{equ}
where $Y^\eps_{x}$ denotes the solution to the second equation
in~\eqref{e:fs_intro} with $X^\eps$ kept fixed at the value~$x$.  If
$\eps $ is small enough that $\Delta t / \eps$ is larger than the
mixing time of~$Y_x^\eps$, the Birkhoff integral in~\eqref{eq:B} is in
fact close to the averaged coefficient in~\eqref{eq:2}, in the sense
that
\begin{equ}
  \label{eq:B}
  F(x) \approx \frac{1}{\Delta t}\int_{n\Delta t}^{(n+1)\Delta t}
  f(x,Y^\eps_{x}(s))ds\;.
\end{equ}
Therefore \eqref{e:scheme1} can also be thought of as an integrator
for the averaged equation~\eqref{eq:2}. In fact, when $\eps$ is small,
one can obtain a good approximation of $F(x)$ using only a fraction of
the macro time step. In particular, we expect that
\begin{equ}
  \label{e:intro_ergodic_approx}
  \frac{1}{\Delta t}\int_{n\Delta t}^{(n+1)\Delta t}
  f(x,Y^\eps_{x}(s))ds \approx \frac{\lambda}{\Delta
    t}\int_{n\Delta t}^{(n+1/\lambda)\Delta t} f(x,Y^\eps_{x}(s))ds =: F_n(x)
\end{equ}
with $\lambda\ge 1$ provided that $\Delta t / (\eps \lambda)$ remains
larger than the mixing time of $Y^\eps_x$. This observation is at the
core of HMM-type methods -- in essence, they amount to replacing
\eqref{e:scheme1} by
\begin{equ}
  \label{e:scheme2}
  x_{n+1} = x_{n} + \Delta t\, F_n(x_n) \;.
\end{equ}
Since the number of computations required to compute the effective
vector field $F_n(x)$ is reduced by a factor~$\lambda$, this is also
the speed-up factor for an HMM-type method.

From the argument above, it is apparent that there is another,
equivalent way to think about HMM-type methods, as was first pointed
out in \cite{fatkullin04} (see
also~\cite{vanden2007hmm,weinan2009general,ariel2012multiscale,%
  ariel2013multiscale}.  Indeed, the integral defining $F_n(x)$
in~\eqref{e:intro_ergodic_approx} can be recast into an integral on
the full interval $[n\Delta t,(n+1)\Delta t]$ by a change of
integration variables, which amount to rescaling the internal clock of
the variables $Y^\eps_x$. In other words, HMM-type methods can also be
thought of as approximating the fast-slow system in~\eqref{e:fs_intro}
by
\begin{equs}
  \label{e:cp_intro}
  \frac{d\Xtilde^\eps}{dt} &= f(\Xtilde^\eps,\Ytilde^\eps) \\
  \frac{d\Ytilde^\eps}{dt} &= \frac{1}{\eps\lambda} g(\Xtilde^\eps,\Ytilde^\eps) \;.
 \end{equs} 
 If $\eps\ll 1$, we can reasonably replace $\eps$ with $\eps \lambda$,
 provided that this product still remains small -- in particular, the
 evolution of the slow variables in~\eqref{e:cp_intro} is still
 captured by the limiting equation~\eqref{eq:2}. Hence HMM-type
 methods are akin to artificial compressibility \cite{chorin67} in
 fluid simulations and Car-Parrinello methods \cite{car85} in
 molecular dynamics.

 The approximations in~\eqref{e:intro_ergodic_approx} or
 \eqref{e:cp_intro} are perfectly reasonable if we are only interested
 in staying faithful to the averaged equation~\eqref{eq:2} -- that is
 to say, HMM-type approximations will have the correct law of large
 numbers (LLN) behavior. However, the fluctuations about that average
 will be enhanced by a factor of $\lambda$. This is quite clear from
 the interpretation \eqref{e:cp_intro}, since in the original model
 \eqref{e:fs_intro}, the local fluctuations about the average are of
 order $\sqrt{\eps}$ and in \eqref{e:cp_intro} they are of order
 $\sqrt{\eps \lambda}$. The large fluctuations about the average
 caused by rare events are similarly inflated by a factor of
 $\lambda$. This can be an issue, for example in metastable fast-slow
 systems, where the large fluctuations about the average determine the
 waiting times for transitions between metastable states. In
 particular we shall see that an HMM-type scheme drastically decreases
 these waiting times due to the enhanced fluctuations.

 In this article we propose a simple modification of HMM which
 corrects the problem of enhanced fluctuations. The key idea is
 to replace the approximation \eqref{e:intro_ergodic_approx} with
\begin{equ}
  \label{e:intro_ergodic_approx2}
  \frac{1}{\Delta t}\int_{n\Delta t}^{(n+1)\Delta t} f(x,Y^\eps_{x}(s))ds
  \approx 
  \sum_{j=1}^\lambda\frac{1}{\Delta
    t}\int_{n\Delta t}^{(n+1/\lambda)\Delta t} 
  f(x,Y^{\eps,j}_{x}(s))ds\;, 
\end{equ}
where each $Y^{\eps,j}_x$ is an independent copy of $Y^\eps_x$. By
comparing \eqref{e:intro_ergodic_approx} with
\eqref{e:intro_ergodic_approx2}, we see that the first approximation
is \emph{essentially} replacing a sum of $\lambda$ weakly correlated
random variables with one random variable, multiplied by
$\lambda$. This introduces correlations that should not be there and
in particular results in enhanced fluctuations. In
\eqref{e:intro_ergodic_approx2}, we instead replace the sum of
$\lambda$ weakly correlated random variables with a sum of $\lambda$
independent random variables. This is a much more reasonable approximation to
make, since these random variables are becoming less and less
correlated as $\eps$ gets smaller. Since the terms appearing on the
right hand side are independent of each other, they can be computed in
parallel. Thus if one has $\lambda$ CPUs available, then the real time
of the computations is identical to HMM. For this reason, we call the
modification the parallelized HMM (PHMM). Note that, in analogy
to~\eqref{e:cp_intro}, one can interpret PHMM as approximating
\eqref{e:fs_intro} by the system
  \begin{equs}
    \label{e:cp_intro_p}
  \frac{d\Xtilde^\eps}{dt} &= \frac{1}{\lambda}\sum_{j=1}^\lambda f(\Xtilde^\eps,\Ytilde^{\eps,j}) \\
  \frac{d\Ytilde^{\eps,j}}{dt} &= \frac{1}{\eps\lambda}
  g(\Xtilde^\eps,\Ytilde^{\eps,j}) \quad \text{for
    $j=1,\dots,\lambda$}\;.
 \end{equs}  
 It is clear that this approximation will be as good as
 \eqref{e:cp_intro} in term of the LLN, but in contrast with
 \eqref{e:cp_intro}, we will show below that it captures the
 fluctuations about the average correctly, both in terms of small
 Gaussian fluctuations and large fluctuations describing rare
 events. A similar observation in the context of numerical
 homogenization was made in~\cite{bal11,bal14}.

 The outline of the remainder of this article is as follows. In
 Section~\ref{s:fs} we recall the averaging principle for stochastic
 fast-slow systems and describe how to characterize the fluctuations
 about this average, including local Gaussian fluctuations and large
 deviation principles. In Section \ref{s:hmm} we recall the HMM-type
 methods. In Section \ref{s:hmm_fluctuations} we show that they lead to
 enhanced fluctuations. In Section \ref{s:phmm} we introduce the PHMM
 modification and in Section \ref{s:phmm_fluctuations} show that this
 approximation yields the correct fluctuations, both in terms of local
 Gaussian fluctuations and large deviations. In Section \ref{s:num} we
 test PHMM for a variety of simple models and conclude in Section
 \ref{s:discussion} with a discussion.

\section{Average and fluctuations in fast-slow systems}
\label{s:fs}

For simplicity we will from here on assume that the fast variables are
stochastic. This assumption is convenient, but not necessary, since
all the averaging and fluctuation properties stated below are known to
hold for large classes of fast-slow systems with deterministically
chaotic fast variables \cite{kifer92, dolgopyat04,
  kelly15a,kelly15b}. The fast-slow systems we investigate are given
by
  \begin{equs}\label{e:fs}
 \frac{dX^\eps}{dt} &= f(X^\eps,Y^\eps) \\
 dY^\eps &= \frac{1}{\eps} g(X^\eps,Y^\eps) dt + \frac{1}{\sqrt{\eps}} \sigma(X^\eps,Y^\eps) dW \;\;,
 \end{equs}
 where $f : \reals^d \times \reals^e \to \reals^d$,
 $g : \reals^d \times \reals^e \to \reals^e$,
 $\sigma: \reals^d \times \reals^e \to \reals^e\times \reals^e$, and
 $W$ is a standard Wiener process in $\reals^e$.  We assume that for
 every $x \in \reals^d$, the Markov process described by the SDE
 \begin{equation}
  dY_x = b(x,Y_x) dt + \sigma(x,Y_x)dW
  \label{eq:1}
\end{equation}
is ergodic, with invariant measure $\mu_x$, and has sufficient mixing
properties. For full details on the necessary mixing properties, see
for instance \cite{fw12}.

In this section we briefly recall the averaging principle for
stochastic fast-slow systems and discuss two results that characterize
the fluctuations about the average, the central limit theorem (CLT)
and the large deviations principle (LDP).

\subsection{Averaging principle}
As $\eps\to 0$, each realization of $X^\eps$, with initial condition
$X^\eps(0) = x$, tends towards a trajectory of a deterministic system
\begin{equ}
  \label{e:fs_avg}
  \frac{d\Xbar}{dt} = F(\Xbar)  \;, \quad \Xbar(0) = x\;,
\end{equ}
where $F(x) = \int f(x,y) \mu_x(dy)$ and $\mu_x$ is the invariant
measure corresponding to the Markov process
$dY_x = b(x,Y_x)dt + \sigma(x,Y_x) dW$. The convergence is in an
almost sure and uniform sense:
\begin{equ}
\lim_{\eps \to 0}\sup_{t \leq T} |X^\eps(t) - \Xbar(t)|  = 0
\end{equ}
for every fixed $T>0$, every choice of initial condition $x$ and
almost surely every initial condition $Y^\eps(0)$ (a.s. with respect
to $\mu_x$) as well as every realization of the Brownian paths driving
the fast variables. Details of this convergence result in the setting
above are given in (for instance) \cite[Chapter 7.2]{fw12}.

\subsection{Small fluctuations -- CLT}
\label{s:fs_clt}

The \emph{small} fluctuations of $X^\eps$ about the averaged system
$\Xbar$ can be understood by characterizing the limiting behavior of
\begin{equ}
Z^\eps : = \frac{X^\eps - \Xbar}{\sqrt{\eps}}\;,
\end{equ}
as $\eps\to 0$. It can be shown that the process $Z^\eps$ converges in
distribution (on the space of continuous functions
$C([0,T]; \reals^d)$ endowed with the sup-norm topology) to a process
$Z$ defined by the SDE
\begin{equ}\label{e:fs_clt}
dZ = B_0(\Xbar) Z dt + \eta (\Xbar) dV \;, \quad  Z(0) = 0\;,
\end{equ} 
Here $\Xbar$ solves the averaged system in~\eqref{e:fs_avg}, $V$ is a
standard Wiener process, $B_0 := B_1 + B_2$ with
\begin{equs}
B_1 (x) &= \int \nabla_x f(x,y) \mu_x(dy )\\
B_2(x) &= \int_0^\infty d\tau \int \mu_x(dy) \nabla_y 
\BE_y \bigg( \ftilde(x,Y_x(\tau))\bigg)  \nabla_x b(x,y)    
\end{equs}
and 
\begin{equ}
\eta(x)\eta^T(x) = \int_0^\infty d\tau\, \BE \ftilde(x, Y_x(0))\otimes 
\ftilde (x,Y_x(\tau) \;,
\end{equ} 
where $\ftilde(x,y) = f(x,y) - F(x)$, $\BE_y$ denotes expectation over
realizations of $Y_x$ with $Y_x(0) = y$, and $\BE$ denotes expectation
over realization of $Y_x$ with $Y_x(0) \sim \mu_x$. We include next a
formal argument deriving this limit, as it will prove useful when
analyzing the multiscale approximation methods. We will replicate the
argument given in \cite{bouchet15}; a more complete and rigorous
argument can be found in \cite[Chapter 7.3]{fw12}.

First, we write a system of equation for the triple
$(\Xbar, Z^\eps, Y^\eps)$ in the following approximated form, which
uses nothing more than Taylor expansions of the original system
in~\eqref{e:fs_intro}:
\begin{equs}
  \frac{d\Xbar}{dt} &= F(\Xbar) \\
  \frac{dZ^\eps}{dt} & = \frac{1}{\sqrt{\eps}} \ftilde(\Xbar, Y^\eps) 
  + \nabla_x f (\Xbar, Y^\eps)  + O(\sqrt{\eps}) \\
  dY^\eps &= \frac{1}{\eps} b(\Xbar,Y^\eps) dt + \frac{1}{\sqrt{\eps}}
  \nabla_x b(\Xbar, Y^\eps) Z^\eps dt + \frac{1}{\sqrt{\eps} }
  \sigma(\Xbar, Y^\eps) dW + O(1)\;.
\end{equs}
We now proceed with a classical perturbation expansion on the
generator of the triple $(\Xbar,Z^\eps,Y^\eps)$. In particular we have
$\CL_\eps = \frac{1}{\eps}\CL_{0} + \frac{1}{\sqrt{\eps}}\CL_1 + \CL_2
+ \dots$ where
\begin{equs}
\CL_0  &=  b(x,y)\cdot\nabla_y + a(x,y) : \nabla_y^2\\
\CL_1 &=  \ftilde(x,y) \cdot\nabla_z + (\nabla_x b(x,y) z) \cdot \nabla_y \\
\CL_2 &= F(x) \cdot \nabla_x + (\nabla_x f(x,y)z) \cdot\nabla_z
\end{equs}
and $a=\sigma\sigma^T$. Let
$u_\eps(x,z,y,t) = \BE_{(x,z,y)} \vphi(\Xbar(t), Z^\eps(t),Y^\eps(t))$
and introduce the ansatz
$u_\eps = u_0 + \sqrt{\eps} u_1 + \eps u_2 + \dots$. By substituting
$u_\eps$ into $\del_t u_\eps = \CL_\eps u_\eps$ and equating powers of
$\eps$ we obtain
\begin{equs}
O(\eps^{-1}) &: \CL_0 u_0 = 0 \\    
O(\eps^{-1/2}) &: \CL_0 u_1 = -\CL_1 u_0 \\
O(\eps^{-1}) &: \del_t u_0 = \CL_2 u_0 + \CL_1 u_1 + \CL_0 u_2\;.
\end{equs}
From the $O(\eps^{-1})$ identity, we obtain $u_0 = u_0 (x,z,t)$,
confirming that the leading order term is independent of $y$. By the
Fredholm alternative, the $O(\eps^{-1/2})$ identity has a solution
$u_1$ which has the Feynman-Kac representation
\begin{equ}
  u_1(x,y,z)  = \int_0^\infty d\tau \, \BE_y \left( \ftilde (x,
    Y_x(\tau) ) \right) 
  \cdot\nabla_z u_0 (x,z)\;,
\end{equ}
where $Y_x$ denotes the Markov process generated by $\CL_0$, i.e. the
solution of~\eqref{eq:1}. Finally, if we average the $O(1)$ identity
against the invariant measure corresponding to $\CL_0$, we obtain
\begin{equs}
\del_t u_0  &= F(x)\nabla_x u_0 + 
\int \mu_x (dy) (\nabla_x f(x,y)z) \cdot\nabla_z u_0  \\  
& + \int \mu_x(dy) \int_0^\infty d\tau \ftilde(x,y) \otimes 
\BE_y\ftilde(x,Y_x(\tau)) : \nabla^2_z u_0 \\
& + \int \mu_x (dy) (\nabla_x b(x,y) z) \int_0^\infty d\tau  \,\nabla_y 
\BE_y\ftilde(x,Y_x(\tau))  \nabla_z u_0 \;.
\end{equs}
Clearly, this is the forward Kolmogorov equation for the Markov
process $(\Xbar, Z)$ defined by
\begin{equs}
\frac{d\Xbar}{dt} &= F(\Xbar) \\
d Z &= B_0(\Xbar) Z dt + \eta (\Xbar ) dV
\end{equs}
with $B_0$ and $\eta$ defined as above. 

\subsection{Large fluctuations -- LDP}
\label{s:fs_ldp}

A large deviation principle (LDP) for the fast-slow system
\eqref{e:fs} quantifies probabilities of $O(1)$ fluctuations of
$X^\eps$ away from the averaged trajectory $\Xbar$. The probability of
such events vanishes exponentially quickly and as a consequence are
not accounted for by the CLT fluctuations, hence an LDP accounts for
the \emph{rare events}.

We say that the slow variables $X^\eps$ satisfy a \emph{large
  deviation principle} (LDP) with action functional $\CS_{[0,T]}$ if
for any set
$\Gamma \subset \{ \gamma \in C([0,T], \reals^d) : \gamma(0) = x \}$
we have
\begin{equs}\label{e:fs_ldp_def}
-\inf_{\gamma \in \mathring{\Gamma}} \CS_{[0,T]}(\gamma) &\leq \liminf_{\eps \to 0} \eps \log \BP \left( X^\eps \in \Gamma \right)\\ &\leq \limsup_{\eps \to 0} \eps \log \BP\left( X^\eps \in \Gamma \right) \leq -\inf_{\gamma \in \bar{\Gamma}} \CS_{[0,T]} (\gamma)\;,
\end{equs}
where $\mathring{\Gamma}$ and $\bar{\Gamma}$ denote the interior and
closure of $\Gamma$ respectively.
 
An LDP also determines many important features of $O(1)$ fluctuations
that occur on large time scales, such as the probability of transition
from one metastable set to another. For example, suppose that $X^\eps$
is known to satisfy an LDP with action functional $\CS_{[0,T]}$. Let
$D$ be an open domain in $\reals^d$ with smooth boundary $\del D$
and let $x^* \in D$ be an asymptotically stable equilibrium for the
averaged system $\dot{\Xbar} = F(\Xbar)$. When $\eps \ll 1$, we expect
that a trajectory of $X^\eps$ that starts in $D$ will tend towards the
equilibrium $x^*$ and exhibit $O(\sqrt{\eps})$ fluctuations about the
equilibrium -- these fluctuations are described by the CLT. On very
large time scales, these small fluctuations have a chance to `pile up'
into an $O(1)$ fluctuation, producing behavior of the trajectory that
would be considered impossible for the averaged system. Such
fluctuations are not accurately described by the CLT and requires the
LDP instead. For example, the asymptotic behaviour of escape time from
the domain $D$,
\begin{equ}
\tau^\eps = \inf \{t > 0 : X^\eps(t) \notin D \}\;,
\end{equ}
can be quantified in terms of the \emph{quasi-potential} defined by
\begin{equ}\label{e:quasi_def}
\CV (x,y) = \inf_{T > 0} \inf_{\gamma(0)=x , \gamma(T) = y} \CS_{[0,T]} (\gamma)
\end{equ}
Under natural conditions, it can be shown that for any $x\in D$
\begin{equ}
\lim_{\eps \to 0} \eps \log \BE_x \tau^\eps = \inf_{y \in \del D} \CV (x^*, y) \;.
\end{equ}
Hence the time it takes to pass from the neighborhood of one
equilibrium to another may be quantified using the LDP. Details on the
escape time of fast-slow systems can be found in \cite[Chapter
7.6]{fw12}.

LDPs for fast-slow systems of the type \eqref{e:fs} are well
understood \cite[Chapter 7.4]{fw12}. First define the Hamiltonian
$\CH : \reals^d \times \reals^d \to \reals$ by
\begin{equ}\label{e:fs_ldp_ham}
\CH (x,\theta) = \lim_{T\to \infty} \frac{1}{T} \log \BE_y \exp \bigg( \theta \cdot \int_0^T f(x,Y_x(s)) ds   \bigg) \;,
\end{equ} 
where $Y_{x}$ denotes the Markov process governed by $dY_x = b(x,Y_x) dt + \sigma(x,Y_x) dW$. Let $\CL : \reals^d \times \reals^d \to \reals$ be the Legendre transform of $\CH$:
\begin{equ}\label{e:fs_ldp_legendre}
\CL (x, \beta) = \sup_{\theta} \left( \theta \cdot \beta - \CH (x, \theta)    \right) \;.
\end{equ}
Then the action functional is given by
\begin{equ}\label{e:action_legendre}
\CS_{[0,T]} (\gamma) = \int_0^T \CL (\gamma(s), \dot{\gamma}(s)) ds\;.
\end{equ} 
It can also be shown that the function
$u(t,x) = \inf_{\gamma(0) = x }\CS_{[0,t]} (\gamma)$ satisfies the
Hamilton-Jacobi equation
\begin{equ}\label{e:fs_ldp_hj}
\del_t u(t,x) = \CH (x ,\nabla u(t,x)) \;.
\end{equ}
Donsker-Varadhan theory tells us that the connection between
Hamilton-Jacobi equations and LDPs is in fact much deeper. Firstly,
\emph{Varadhan's Lemma} states that if a process $X^\eps$ is known to
satisfy an LDP with some associated Hamiltonian $\CH$, then for any
$\vphi : \reals^d \to \reals$ we have the generalized Laplace
method-type result
\begin{equ}\label{e:fs_ldp_varadhan}
\lim_{\eps \to 0 } \eps \log \BE_x \exp \left( \eps^{-1} \vphi (X^\eps (t))   \right) = S_t \vphi (x)
\end{equ}
where $S_t$ is the semigroup associated with the Hamilton-Jacobi
equation $\del_t u = \CH (x,\nabla u)$. Conversely, if it is known
that \eqref{e:fs_ldp_varadhan} holds for all $(x,t)$ and a suitable
class of $\vphi$, then the inverse Varadhan's lemma states that
$X^\eps$ satisfies an LDP with action functional given by
\eqref{e:fs_ldp_legendre}, \eqref{e:action_legendre}. Hence we can use
\eqref{e:fs_ldp_varadhan} to determine the action functional for a
given process.

In the next few sections, we will exploit both sides of Varadhan's
lemma when investigating the large fluctuations of the HMM and related
schemes. More complete discussions on Varadhan's Lemma can be found in
\cite[Chapters 4.3, 4.4]{dz09}.

 \section{HMM for fast-slow systems}
\label{s:hmm}
When applied to the stochastic fast-slow system \eqref{e:fs}, HMM-type
schemes rely on the fact that the slow $X^\eps$ variables, and the
coefficients that govern them, converge to a set of reduced variables
as $\eps$ tends to zero. We will describe a simplest version of the
method below, which is more convenient to deal with mathematically. 

Before proceeding, we digress briefly on notation. When referring to
continuous time variables we will always use upper case symbols
($X^\eps,Y^\eps$ etc) and when referring to discrete time
approximations we will always use lower case symbols ($x^\eps_n$,
$y^\eps_n$ etc). We will also encounter continuous time variables
whose definition depends on the integer $n$ for which we have
$t \in [n\Delta t, (n+1)\Delta t)$. We will see below that such
continuous time variables are used to define discrete time
approximations. In this situation we will use upper case symbols with
a subscript $n$ (eg. $X^\eps_n$).

Let us now describe a `high-level' version of HMM.  Fix a step size
$\Delta t$ and define the intervals
$I_{n, \Delta t}: = [n\Delta t, (n+1)\Delta t)$. On each interval
$I_{n,\Delta t}$ we update $x^\eps_n \approx X^\eps(n\Delta t)$ to
$x^\eps_{n+1} \approx X^\eps((n+1)\Delta t)$ via an iteration of the
following two steps:
 \begin{enumerate}
 \item (Micro step) Integrate the fast variables over the interval
   $I_{n,\Delta t}$, with the slow variable frozen at
   $X^\eps = x^\eps_n$. That is, the fast variables are approximated
   by
\begin{equ}\label{e:micro}
Y^\eps_{n}(t) = Y^\eps_n(n\Delta t) + \frac{1}{\eps} 
\int_{n\Delta t}^t g(x^\eps_n,Y_{n}^\eps (s))ds 
+  {\frac{1}{\sqrt{\eps}}}\int_{n\Delta t}^t \sigma (x^\eps_n, Y^{\eps}_n(s)) dW(s) 
\end{equ}
for $n\Delta t \leq t \leq (n+ 1/\lambda)\Delta t $ with some
$\lambda \geq 1$ (that is, we do not necessarily integrate the
$Y_n^\eps$ variables over the whole time window). Due to the
ergodicity of $Y_x$, the initialization of $Y^\eps_n$ is not crucial
to the performance of the algorithm. It is however convenient to use
$Y^\eps_{n+1}(0) = Y^\eps_n((n+ 1/\lambda)\Delta t)$, since this
reinitialization leads to the interpretation of the HMM scheme given
in~\eqref{e:hmm_cp} below.

\item (Macro step) Use the time series from the micro step to update
  $x^\eps_n$ to $x^\eps_{n+1}$ via
\begin{equ}
  \label{e:macro}
  x^{\eps}_{n+1} = x^\eps_n + \lambda \int_{n\Delta
    t}^{(n+1/\lambda)\Delta t} f(x^\eps_n,Y^\eps_n(s)) ds\;.
\end{equ}
Note that we do not require $Y^\eps_n$ over the whole $\Delta t$ time
step, but only a fraction of the step large enough for $Y^\eps_n$ to
mix. Indeed, if $\eps$ is small enough, we have the approximate equality
\begin{equ}
\frac{\lambda}{ \Delta t}\int_{n\Delta t}^{(n+1/\lambda)\Delta
  t} 
f(x^\eps_n,Y^\eps_n(s)) ds \approx
\frac{1}{\Delta t}\int_{n\Delta t}^{(n+1)\Delta t} f(x^\eps_n,Y^\eps_n(s))
ds 
\end{equ}
since both sides are close the the ergodic mean
$\int f(x^\eps_n , y) d\mu_{x^\eps_n}(y )$. 
\end{enumerate} 
Clearly, the efficiency of the methods comes from the fact that we do
not need to compute the fast variables on the whole time interval
$I_{n,\Delta t}$ but only a $1/\lambda$ fraction of it. Hence
$\lambda$ should be considered the \emph{speed-up factor} of HMM.

As already stated, the algorithm above  is a high-level version,
in that one must do further approximations to make the method
implementable. For example, one typically must specify some
approximation scheme to integrate \eqref{e:micro}, for instance with
Euler-Maruyama we compute the time series by
\begin{equ}\label{e:micro_euler}
y^\eps_{n,m+1} = y^\eps_{n,m} + \frac{\delta t}{\eps}
g(x^\eps_n,y^\eps_{n,m}) 
+ \sqrt{\frac{\delta t}{\eps}} \sigma (x^\eps_n, y^\eps_{n,m}) \xi_{n,m} \;,
\end{equ}
where $0 \leq m \leq M$ is the index within the micro step,
$\xi_{n,m}$ are i.i.d. standard Gaussians and the micro-scale step
size $\delta t$ is much smaller than the macro-scale step size
$\Delta t$. In the macro step, we would similarly have
\begin{equ}\label{e:macro_euler}
x^\eps_{n+1} = x^\eps_n + \Delta t  \, F_n(x_n^\eps) 
\end{equ}
where $F_n(x) = \frac{1}{M} \sum_{m=1}^M f(x, y^\eps_{n,m})$ and
$ M = \Delta t / (\delta t \lambda)$.

The following observation, which is taken from \cite{fatkullin04},
will allow us to easily describe the average and fluctuations of the
above method. On each interval $I_{n,\Delta t}$, the high-level HMM
scheme described above is equivalently given by
$x^\eps_{n+1} = X^\eps_n((n+1)\Delta t)$, where $X^\eps_n$ solves the
system
\begin{equs}\label{e:hmm_cp}
\frac{dX^\eps_n}{dt} &= f(x^\eps_n,\Ytilde^\eps_n) \\
d\Ytilde^\eps_n &= \frac{1}{\eps \lambda} b(x^\eps_n,\Ytilde^\eps_n) dt + \frac{1}{\sqrt{\eps \lambda}}\sigma (x^\eps_n,\Ytilde^\eps_n) dB\;,
\end{equs}
defined on the interval $n\Delta t \leq t \leq (n+1)\Delta t$, with
the initial condition $X^\eps_n(n\Delta t) = x^\eps_n$. This can be
checked by a simple rescaling of time. It is clear that the efficiency
of HMM essentially comes from saying that the fast-slow system is not
drastically changed if one replaces $\eps$ with the slightly larger,
but still very small $\eps \lambda$. 

\section{Average and fluctuations in HMM methods}
\label{s:hmm_fluctuations}

In this section we investigate whether the limit theorems discussed
in Section \ref{s:fs}, i.e. the averaging principle, the CLT
fluctuations and the LDP fluctuations, are also valid in the HMM
approximation a fast-slow system. We will see that the averaging
principle is the only property that holds, and that both types of
fluctuations are \emph{inflated} by the HMM method.

\subsection{Averaging}\label{s:hmm_avg}
By construction, HMM-type schemes capture the correct averaging
principle. More precisely, if we take $\eps \to 0$ then the sequence
$x^\eps_n$ converges to some $\xbar_n$, where $\xbar_n$ is a numerical
approximation of the true averaged system $\Xbar$. If this numerical
approximation is well-posed, the limits $\eps \to 0$ and
$\Delta t \to 0$ commute with one another. Hence the HMM approximation
$x^\eps_n$ is consistent, in that it features approximately the same
averaging behavior as the original fast-slow system.

We will argue the claim by induction. Suppose that for some $n\geq 0$ we know that $\lim_{\eps \to 0} x^\eps_n = \xbar_n$  (the $n=0$ claim is trivial, since they are both simply the initial condition). Then, using the representation \eqref{e:hmm_cp} we know that $x^\eps_{n+1} = X^\eps_n((n+1)\Delta t)$ where $X^\eps_n(n\Delta t) = x^\eps_n$. Since \eqref{e:hmm_cp} is a fast-slow system of the form \eqref{e:fs} we can apply the averaging principle from Section \ref{s:fs}. In particular it follows that $X^\eps_n \to \Xbar_n$ uniformly (and almost surely) on $I_{n,\Delta t}$, where $\Xbar_n$ satisfies the averaged ODE
\begin{equ}
\frac{d\Xbar_n}{dt} = \int f(\xbar_n , y) \mu_{ \xbar_n} (dy ) = F(\xbar_n)\;.
\end{equ}  
Since the right hand side is a constant, it follows that $x^\eps_{n+1} \to \xbar_{n+1}$ as $\eps \to 0$, where
\begin{equ}
\xbar_{n+1} = \xbar_n + F( \xbar_n) \Delta t \;.
\end{equ} 
This is nothing more than the Euler approximation of the true averaged variables $\Xbar$, which completes the induction and hence the claim. 
\par
Introducing an integrator in to the micro-step will make things more complicated, as the invariant measures appearing will be those of the discretized fast variables. In \cite{mattingly02} it is shown that discretizations of SDEs often do not possess the ergodic properties of the original system. For those situations where no such issues arise, rigorous arguments concerning this scenario, including rates of convergence for the schemes, are given in \cite{liu05}.

\subsection{Small fluctuations}\label{s:hmm_clt}
For HMM-type methods, the CLT fluctuations about the average become inflated by a factor of $\sqrt{\lambda}$. That is, if we define
\begin{equ}
z^\eps_{n+1} = \frac{x^\eps_{n+1} - \xbar_{n+1}}{\sqrt{\eps}} \;,
\end{equ}
then as $\eps\to 0$, the fluctuations described by $z^\eps_{n+1}$ are not consistent with \eqref{e:fs_clt}, but rather with the SDE 
\begin{equ}\label{e:hmm_clt_z}
dZ = B(\Xbar) Z dt + \sqrt{\lambda} \eta(\Xbar) dV\;, \quad Z(0) = 0
\end{equ} 
where $\Xbar$ satisfies the correct averaged system. 
\par
As above, by consistency we mean that when we take $\eps\to 0$, the sequence $\{z^\eps_n\}_{n\geq 0}$ converges to some well-posed discretization of the SDE \eqref{e:hmm_clt_z}. Since $Z(0) = 0$, it is easy to see that the solution to this equation is simply $\sqrt{\lambda}$ times the solution of \eqref{e:fs_clt}. Hence the fluctuations of the HMM-type scheme are inflated by a factor of $\sqrt{\lambda}$.  
\par
It is convenient to look instead at the rescaled fluctuations 
\begin{equ}
\zhat^\eps_n = z^\eps_n / \sqrt{\lambda} =  \frac{x^\eps_n - \xbar_n}{\sqrt{\eps \lambda}}\;,
\end{equ}
since this allows us to reproduce the argument from Section \ref{s:fs_clt}, with $\eps ' = \eps \lambda$ playing the role of $\eps$. We will again argue by induction, assuming for some $n\geq 0$ that $\zhat^\eps_n \to \zhat_n$ as $\eps \to 0$ (the $n=0$ case is trivial). 
\par
The rescaled fluctuations are given by $\zhat^\eps_{n+1} = Z^\eps_n((n+1)\Delta t)$ where $Z^\eps_n (t) = (X^\eps_n(t) - \Xbar_n(t)) / \sqrt{\eps \lambda}$ and $X^\eps_n(t)$ is governed by the system \eqref{e:hmm_cp} with initial condition $X^\eps_n(n\Delta t) = x^\eps_n$ and $\Xbar_n$ satisfies
\begin{equ}
\frac{d\Xbar_n}{dt} = F(\xbar_n) 
\end{equ}
with initial condition $\Xbar_n (n\Delta t) = \xbar_n$. We can then obtain the reduced equations for the pair $(X_n^\eps, Z^\eps_n)$ by arguing exactly as in Section \ref{s:fs}. Indeed, the triple $(\Xbar_n, Z^\eps_n, \Ytilde^\eps_n)$ is governed by the system 
\begin{equs}
\frac{d\Xbar_n}{dt} &= F(\xbar_n)\\
\frac{d \Zhat^\eps_n}{dt} & = \frac{1}{\sqrt{\eps \lambda }} \ftilde (\xbar_n, \Ytilde^\eps_n) + \nabla_x f (\xbar_n , \Ytilde^\eps_n) \zhat_n + O(\sqrt{\eps \lambda })	   \\
d\Ytilde^\eps_n & = \frac{1}{\eps \lambda} b(\xbar_n , \Ytilde^\eps_n) dt + \frac{1}{\sqrt{\eps \lambda }} \nabla_x b(\xbar_n , \Ytilde^\eps_n) \zhat_n dt + \frac{1}{\sqrt{\eps \lambda }}\sigma (\xbar^\eps_n, \Ytilde^\eps_n) dW + O(1)
\end{equs}
From here on we can carry out the calculation precisely as in Section \ref{s:fs_clt}, with the added convenience of the vector fields no longer depending on $x$ as a variable. In doing so we obtain $\Zhat^\eps_n \to \Zhat_n$ (in distribution) as $\eps \to 0$, where
\begin{equ}
d\Zhat_n =  B_0(\xbar_n) \zhat_n dt + \eta(\xbar_n) dV\;,
\end{equ}   
with the initial condition defined recursively by $\Zhat_n (n\Delta t) =\zhat_n$. Using the fact that $\zhat_{n+1} = \Zhat_n((n+1)\Delta t)$, we obtain
\begin{equ}
\zhat_{n+1} = \zhat_n + B_0 (\xbar_n) \zhat_n \Delta t + \eta(\xbar_n) \sqrt{\Delta t} \xi_n
\end{equ} 
where $\xi_n$ are iid standard Gaussians. Hence we obtain the
Euler-Maruyama scheme for the correct CLT \eqref{e:fs_clt}. However,
since $\zhat^\eps_n$ describes the rescaled fluctuations, we see that
the true fluctuations $z^\eps_n$ of HMM are consistent with the
inflated \eqref{e:hmm_clt_z}.

\subsection{Large fluctuations}
\label{s:hmm_ldp}

As with the CLT, the LDP of the HMM scheme is not consistent with the true LDP of the fast-slow system, but rather a rescaled version of the true LDP. In particular, define $u_{\lambda,\Delta t}$ by   
\begin{equ}
u_{\lambda , \Delta t}(t,x) = \lim_{\eps \to 0} \eps \log \BE_{x} \exp \bigg(\frac{1}{\eps} \vphi( x^\eps_{n+1} ) \bigg)
\end{equ}
for $t \in I_{n,\Delta t}$. If the $O(1)$ fluctuations of HMM were consistent with those of the fast-slow system, we would expect $u_{\lambda,\Delta}$ to converge to the solution of \eqref{e:fs_ldp_hj} as $\Delta t \to 0$. Instead, we find that  as $\Delta t\to 0$,  $u_{\lambda,\Delta t}(t,x)$ converges to the solution to the Hamilton-Jacobi equation
\begin{equ}\label{e:hmm_ldp}
\del_t u_\lambda  = \frac{1}{\lambda} \CH (x ,\lambda \nabla u_\lambda)\;\quad\; \quad u_\lambda(0,x) = \vphi(x) \;.
\end{equ}
In light of the discussion in Section \ref{s:fs_ldp}, the reverse Varadhan lemma suggests that the HMM scheme is consistent with the wrong LDP. Before proving this claim, we first discuss some implications.  
\par
{{The rescaled Hamilton-Jacobi equation implies that the action functional for HMM will be a rescaled version of that for the true fast-slow system. 
%
%
Indeed, it is easy to see that the Langrangian corresponding to HMM simplifies to 
\begin{equ}
\widehat{\CL}(x,\beta) := \sup_{\theta} \left( \theta \cdot \beta - \frac{1}{\lambda} \CH(x,\lambda\theta) \right) = \frac{1}{\lambda} \CL(x,\beta)\;,
\end{equ}
where $\CL$ is the Lagrangian for the true fast-slow system. Thus, the action of the HMM approximation is given by $\widehat{\CS}_{[0,T]} = \lambda^{-1} \CS_{[0,T]}$ where $\CS$ is the action of the true fast-slow system. \par
In particular, it follows immediately from the definition that the HMM approximation has quasi-potential $\widehat{\CV}(x,y) = \lambda^{-1} \CV(x,y)$, where $\CV$ is the true quasi-potential. As a consequence, the escape times for the HMM scheme will be
drastically faster than those of the fast-slow system. In the
terminology of Section \ref{s:fs_ldp}, if we let
$\tau^{\eps,\Delta t}$ be the escape time for the HMM scheme then for
$\eps, \Delta t \ll 1$ we expect
\begin{equ}\label{e:hmm_ldp_escape}
  \BE \tau^{\eps,\Delta t} \asymp\exp\Big( \frac{1}{\eps \lambda}\CV
  (x^*, \del D) \Big) \;.
\end{equ} 
where $\asymp$ log-asymptotic equality.  Thus, the log-expected
escape times are decreasing proportionally with $\lambda$. On the other hand, since the HMM action is a multiple of the true action, the minimizers will be unchanged by the HMM approximation. Hence the large deviation transition pathways will be unchanged by the HMM approximation. }}

%
%
%
%
To justify the claim for $u_{\lambda,\Delta t}$ \eqref{e:hmm_ldp}, we first introduce
some notation. Let $S_{t}^{(\alpha)}$ be the semigroup associated with
the Hamilton-Jacobi equation
\begin{equ}\label{e:hmm_hj_alpha}
\del_t v(t,x) = \CH (\alpha, \nabla v(t,x))\;,
\end{equ}
notice that this is the same as the true Hamilton-Jacobi equation
\eqref{e:fs_ldp_hj} but with the first argument of the Hamiltonian now
frozen as a parameter $\alpha$. The necessity of the parameter
$\alpha$ is due to the fact that in the system for
$(X^\eps_n, Y^\eps_n)$, the $x$ variable in the fast process is frozen
to its value at the left end point of the interval, and hence is
treated as a parameter on each interval.
We also introduce the operator
$S_{t} \psi (x) = S^{(\alpha)}_{t} \psi (x) |_{\alpha = x}$ and also
$S_{\lambda, t} = \lambda^{-1} S_{ t} (\lambda \cdot)$. In this
notation, it is simple to show that
\begin{equ}\label{e:hmm_ldp_Qn}
\BE_x \exp \left( \eps^{-1} \vphi (x^\eps_n)  \right) \asymp \exp\left(\eps^{-1} (S_{\lambda, \Delta t})^n \vphi(x) \right)\;.
\end{equ}
We will verify \eqref{e:hmm_ldp_Qn} by induction, starting with the
$n=1$ case. Since, on the interval $I_{0,\Delta t}$, the pair
$(X^\eps_0, \Ytilde^\eps_0)$ is a fast-slow system of the form
\eqref{e:fs} with $\eps$ replaced by $\eps \lambda$, it follows from
Section \ref{s:fs_ldp} that $X^\eps_0$ satisfies an LDP with action
functional derived from the Hamiltonian-Jacobi equation
\eqref{e:hmm_hj_alpha}, with the parameter $\alpha$ set to the value
of $X^\eps_0$ at the left endpoint, which is $X^\eps_0(0) = x$. Hence,
it follows from Varadhan's lemma that for any suitable
$\psi : \reals^d \to \reals$
\begin{equ}
\BE_x \exp \left((\eps \lambda)^{-1} \psi (X^\eps_0(\Delta t))\right) \asymp \exp \left( (\eps \lambda)^{-1} S_{\Delta t}^{(\alpha)} \psi(x)|_{\alpha = x}  \right) \;.
\end{equ}
Hence, since $x_1^\eps = X^\eps_0(\Delta t)$ with $X^\eps_0(0) = x$, we have
\begin{equs}
\BE_x \exp \left(\eps^{-1} \vphi (x^\eps_1) \right) &= \BE_x \exp \left((\eps\lambda)^{-1} \lambda \vphi (X^\eps_1(\Delta t))\right)\\ &\asymp \exp \left((\eps \lambda)^{-1} S^{(\alpha)}_{\Delta t}(\lambda \vphi)(x)|_{\alpha = x} \right) = \exp \left(\eps^{-1} S_{\lambda, \Delta t}\vphi(x) \right)
\end{equs}
as claimed. Now, suppose \eqref{e:hmm_ldp_Qn} holds for all $k$ with $n \geq k \geq 1$, then 
\begin{equ}\label{e:hmm_lpd_Q1a}
\BE_x \exp \left( \eps^{-1} \vphi (x^\eps_{n+1})  \right)  = \BE_x \BE_{x_1^\eps} \exp \left( \eps^{-1} \vphi (x^\eps_{n+1})  \right) \;.
\end{equ} 
By the inductive hypothesis, we have that
\begin{equ}\label{e:hmm_lpd_Q1b}
\BE_{x_1^\eps} \exp \left( \eps^{-1} \vphi (x^\eps_{n+1})  \right) \asymp \exp\left(\eps^{-1} (S_{\lambda,\Delta t})^n \vphi(x^\eps_1) \right)\;.
\end{equ}
Applying \eqref{e:hmm_lpd_Q1b} under the expectation in
\eqref{e:hmm_lpd_Q1a} (see Remark~\ref{rmk:under_exp}) we see that
\begin{equ}
\BE_x \exp \left( \eps^{-1} \vphi (x^\eps_{n+1})  \right) = \BE_x \BE_{x^\eps_1} \exp \left( \eps^{-1} \vphi (x^\eps_{n+1})  \right)  \asymp \BE_x \exp\left(\eps^{-1} (S_{\lambda,\Delta t})^n \vphi(x^\eps_1) \right)\;.
\end{equ}
Now applying the inductive hypothesis with $n=1$ and $\psi (\cdot) = (S_{\lambda,\Delta t})^n \vphi(\cdot) $
\begin{equ}
\BE_x \exp\left(\eps^{-1} (S_{\lambda,\Delta t})^n \vphi(x^\eps_1) \right) \asymp \exp\left(\eps^{-1} S_{\lambda,\Delta t} (S_{\lambda,\Delta t})^n \vphi(x) \right)\;,
\end{equ}
which completes the induction.
\par
By definition, we therefore have $u_{\lambda,\Delta t}(t,x) = (S_{\lambda,\Delta t})^n \vphi (x)$ when $t \in I_{n,\Delta t}$. All that remains is to argue that $u_{\lambda,\Delta t}$ converges to the solution of \eqref{e:hmm_ldp} as $\Delta t \to 0$. But this can be seen from the expansion of the semigroup 
\begin{equs}
\frac{u_{\lambda,\Delta t}(t + \Delta t,x) - u_{\lambda,\Delta t}(t,x)}{\Delta t} &= 
\frac{(S_{\lambda,\Delta t})^{n+1}\vphi (x)  - (S_{\lambda,\Delta t})^{n}\vphi (x)}{\Delta t}\label{e:hmm_ldp_expansion} \\  &= \frac{S_{\lambda,\Delta t} (S_\Delta t)^n \vphi (x) - (S_{\lambda,\Delta t})^n \vphi (x)}{\Delta t}\\ &= \lambda^{-1} \CH (\alpha , \lambda \nabla (S_{\lambda,\Delta t} )^n \vphi(x) )|_{\alpha = x} + O(\Delta t) \\ 
& = \lambda^{-1} \CH (x , \lambda \nabla u_{\lambda,\Delta t}(t,x))) + O(\Delta t)
\end{equs}
which yields the desired limiting equation. 

 \begin{rmk}
\label{rmk:under_exp}
Regarding the operation of taking the log-asymptotic result inside the
expectation, one can find such calculations done rigorously in (for
instance) \cite[Lemma 4.3]{fw12}.
 \end{rmk}

\begin{rmk}
\label{rmk:rescaling}
From the discussion above, it appears that the mean transition time
can be estimated from HMM upon exponential rescaling, see
\eqref{e:hmm_ldp_escape}. This is true, but only at the level of the
(rough) log-asymptotic estimate of this time. How to rescale the
prefactor is by no means obvious. As we will see below PHMM avoids
this issue altogether since it does not necessitate any rescaling.
 \end{rmk}

 \section{Parallelized HMM}
\label{s:phmm}
 
 There is a simple variant of the above HMM-type scheme which captures
 the correct average behavior and  fluctuations, both at the level of the CLT
 and LDP. In a usual HMM type method, the key approximation is given
 by
\begin{equ}\label{e:hmm_approx}
\int_{n\Delta t}^{(n+1)\Delta t} f(x^\eps_n , Y^\eps_{n}(s))ds \approx \lambda \int_{n\Delta t}^{(n+1/\lambda)\Delta t} f(x^\eps_n , Y^\eps_{n}(s))ds\;,
\end{equ}
which only requires computation of the fast variables on the interval $[n\Delta t , (n+1/\lambda)\Delta t]$. This approximation is effective at replicating averages, but poor at replicating fluctuations. Indeed,  for each $j$, the time series $Y^\eps_{n}$ on the interval $[(n+j/\lambda)\Delta t ,(n+(j+1)/\lambda)\Delta t]$ is replaced with an identical copy of the time series from the interval $[n \Delta t ,(n+1/\lambda)\Delta t]$. This introduces strong correlations between random variables that should be essentially independent. Parallelized HMM avoids this issue by employing the approximation
\begin{equ}
\int_{n\Delta t}^{(n+1)\Delta t} f(x^\eps_n , Y^\eps_{n}(s))ds \approx \sum_{j=1}^\lambda \int_{n\Delta t}^{(n+1/\lambda)\Delta t} f(x^\eps_n , Y^{\eps,j}_n(s))ds\;,
\end{equ}
where $Y^{\eps,j}_n$ are for each $j$ independent copies of the time series computed in \eqref{e:hmm_approx}. Due to their independence, each copy of the fast variables can be computed in parallel, hence we refer to the method as parallel HMM (PHMM). The method is summarized below. 
 
  \begin{enumerate}
  \item (Micro step) On the interval $I_{n,\Delta t}$, simulate
    $\lambda$ independent copies of the of the fast-variables, each
    copy simulated precisely as in the usual HMM. That is, let
\begin{equ}\label{e:microp}
Y^{\eps,j}_n = Y^{\eps,j}_n(n\Delta t) + \frac{1}{\eps} \int_{n\Delta t}^t g(x^\eps_n,Y^{\eps,j}_n(s))ds +  {\frac{1}{\sqrt{\eps}}}\int_{n\Delta t}^t \sigma (x^\eps_n, Y^{\eps,j}_n(s)) dW_j(s) 
\end{equ}
for $j=1,\dots,\lambda$ with $W_j$ independent Brownian motions. As
with ordinary HMM, we will not require the time series of the whole
interval $I_{n,\Delta t}$ but only over the subset
$[n\Delta t, (n + 1/\lambda)\Delta t )$.
\item (Macro step) Use the time series from the micro step to update
  $x^\eps_n$ to $x^\eps_{n+1}$ by
\begin{equ}\label{e:macrop}
x^\eps_{n+1} = x^\eps_n + \sum_{j=1}^\lambda \int_{n\Delta t}^{(n+1/\lambda)\Delta t} f(x^\eps_n,Y^{\eps,j}_n(s)) ds\;. 
\end{equ}
\end{enumerate} 

As with the HMM-type method, it will be convenient to write PHMM as a
fast-slow system (when restricted to an interval $I_{n,\Delta
  t}$).
Akin to \eqref{e:hmm_cp}, it is easy to verify that the parallel HMM
scheme is described by the system
\begin{equs}\label{e:phmm_cp}
\frac{dX^\eps_{n}}{dt} &= \frac{1}{\lambda}\sum_{j=1}^\lambda f(x^\eps_n,\Ytilde^{\eps,j}_n) \\
d\Ytilde^\eps_{n,j} &= \frac{1}{\eps \lambda} b(x^\eps_n,\Ytilde^{\eps,j}_n) dt + \frac{1}{\sqrt{\eps \lambda}}\sigma (x^\eps_n,\Ytilde^{\eps,j}_n) dW_j\;,
\end{equs}
for $j=1,\dots,\lambda$ with the initial condition $X^\eps_{n}(n\Delta t) = x^\eps_n$. 

\section{Average and fluctuations in parallelized HMM}
\label{s:phmm_fluctuations}

In this section we check that the averaged behavior and the
fluctuations in the PHMM method are consistent with those in the original fast
slow system.
\subsection{Averaging}
Proceeding exactly as in Section \ref{s:hmm_avg}, it follows that as $\eps \to 0$ the PHMM scheme $x^\eps_{n+1}$ converges to $\xbar_{n+1} = {\Xbar}_{n} ((n+1)\Delta t) $ where
\begin{equ}\label{e:phmm_avg}
\frac{d\Xbar_{n}}{d t} = \frac{1}{\lambda}\sum_{j=1}^\lambda F (\xbar_n) =  F(\xbar_n)
\end{equ}
with initial condition $\Xbar_{n} (n \Delta t) = \xbar_n$. Hence, we
are in the exact same situation as with ordinary HMM, so the averaged
behavior is consistent with that of the original fast slow system.

\subsection{Small fluctuations}\label{s:phmm_clt}
We now show that the fluctuations
\begin{equ}
z^\eps_{n} = \frac{x^\eps_n - \xbar_n}{\sqrt{\eps}}
\end{equ}
are consistent with the correct CLT fluctuations, described by \eqref{e:fs_clt}. As in Section \ref{s:hmm_clt}, we instead look at the rescaled fluctuations
\begin{equ}
\zhat^\eps_{n} = \frac{x^\eps_n - \xbar_n}{\sqrt{\eps \lambda }}\;.
\end{equ}
In particular we will show that these rescaled fluctuations are consistent with 
\begin{equ}\label{e:phmm_clt_a}
d\Zhat = B_0(\Xbar)\Zhat dt + \lambda^{-1/2} \eta(\Xbar) dV\;.
\end{equ}
The claim for $z^{\eps}$ will follows immediately from the claim for $\zhat^\eps$. 
\par
We have that $\zhat^\eps_{n+1} = \Zhat^\eps_{n} ((n+1)\Delta t)$ where
\begin{equ}
\Zhat^\eps_{n}(t)=  \frac{X^\eps_{n}(t) - \Xbar_{n}(t)}{\sqrt{\eps \lambda}} 
\end{equ} 
with $X^\eps_{n}$ given by the system \eqref{e:phmm_cp} and $\Xbar_{n}$ given by the averaged equation \eqref{e:phmm_avg}. As in Section \ref{s:hmm_clt}, we derive a system for the triple $(\Xbar_{n}, \Zhat^\eps_{n}, \Ytilde^\eps_{n})$, where now the fast process has $\lambda$ independent components $\Ytilde^\eps_{n} = (\Ytilde^{\eps,1}_{n}, \dots,\Ytilde^{\eps,\lambda}_{n})$:
\begin{equs}\label{e:phmm_clt_triple}
\frac{d\Xbar_{n}}{dt} &= F(\xbar_n)\\
\frac{d \Zhat^\eps_{n}}{dt} & = \frac{1}{\sqrt{\eps \lambda }} \frac{1}{\lambda}\sum_{j=1}^\lambda \ftilde (\xbar_n, \Ytilde^{\eps,j}_{n}) + \frac{1}{\lambda}\sum_{j=1}^\lambda \nabla_x f (\xbar_n , \Ytilde^{\eps,j}_{n}) \zhat_n + O(\sqrt{\eps \lambda })	   \\
d\Ytilde^{\eps,j}_{n} & = \frac{1}{\eps \lambda} b(\xbar_n , \Ytilde^{\eps,j}_{n}) dt + \frac{1}{\sqrt{\eps \lambda }} \nabla_x b(\xbar_n , \Ytilde^{\eps,j}_{n}) \zhat_n dt + \frac{1}{\sqrt{\eps \lambda }}\sigma (\xbar_n, \Ytilde^{\eps,j}_{n}) dW_j + O(1)\;.
\end{equs}
 With a modicum added difficulty, we can now argue as in Section \ref{s:fs_clt} with $\eps ' = \eps \lambda$ playing the role of $\eps$. The invariant measure $\mu_x^\lambda (dy)$ associated with the generator of $Y^\eps_n$ is now the product measure
\begin{equ}
\mu_x^\lambda (dy_1,\dots,dy_\lambda) = \mu_x(dy_1) \dots \mu_x (dy_\lambda)
\end{equ}
where $\mu_x$ is the invariant measure associated with $\CL_0$ from Section \ref{s:fs_clt}. This product structure simplifies the seemingly complicated expressions arising in the perturbation expansion of \eqref{e:phmm_clt_triple}. In the setting of Section \ref{s:fs_clt} we have that $u_0 = u_0 (x,z,t)$ and 
\begin{equs}\label{e:phmm_u1}
u_1 (x,z,y,t)  = (-\CL_0^{(1)} - \dots - \CL_0^{(\lambda)})^{-1} \CL_1 u_0 (x,z,y,t)\;,
\end{equs}
where $\CL^{(j)}_0 = b(\xbar_n , y_j) \nabla_{y_j} + \frac{1}{2} \sigma \sigma^T (\xbar_n , y_j ) : \nabla_{y_j}^2$

Since 
\begin{equ}
\CL_1 u_0 (x,z,y,t) = \frac{1}{\lambda}\sum_{j=1}^\lambda 
\ftilde (\xbar_n,y_j) \cdot \nabla_z u_0 (x,z,t)\;,
\end{equ}
the Feynman-Kac representation of \eqref{e:phmm_u1} yields
\begin{equ}
u_1(x,z,y,t) = \frac{1}{\lambda} \sum_{j=1}^{\lambda} 
\int_0^\infty d\tau \BE_{y_j} \ftilde (\xbar_n, Y_{\xbar_n,j}(\tau)) 
\cdot \nabla_z u_0(x,z,t)\;.
\end{equ}
The equation for $u_0$ is now given by
\begin{equs}
  \label{e:phmm_u0_big}
  \del_t u_0  &= F(\xbar_n)\nabla_x u_0 + \int \mu_{\xbar_n} (dy_1)\dots \mu_{\xbar_n}( dy_\lambda) (\frac{1}{\lambda}\sum_{j=1}^\lambda \nabla_x f(\xbar_n,y_j)\zhat_n) \nabla_z u_0  \\  
&\quad  + \int \mu_{\xbar_n} (dy_1)\dots \mu_{\xbar_n}( dy_\lambda) \\
& \qquad \times \Biggl( \int_0^\infty d\tau \left(\frac{1}{\lambda} \sum_{j=1}^\lambda \ftilde(\xbar_n,y_j) \right) \otimes \left( \frac{1}{\lambda}\sum_{k=1}^\lambda \BE_y\ftilde(\xbar_n,Y^k_{\xbar_n}(\tau))\right) : \nabla^2_z u_0 \\
& \qquad + \sum_{j=1}^\lambda (\nabla_x b(\xbar_n,y_j) \zhat_n) \int_0^\infty d\tau \, \nabla_{y_j} \frac{1}{\lambda}\sum_{k=1}^\lambda\BE_{y_k}\ftilde(\xbar_n,Y^k_{\xbar_n}(\tau))  \nabla_z u_0 \Biggr)\;.
\end{equs}
By expanding the product measure, the second term on the right hand side of \eqref{e:phmm_u0_big} becomes 
\begin{equ}
  \begin{aligned}
    & \frac{1}{\lambda}\sum_{j=1}^\lambda \int \mu_{\xbar_n}( dy_j)
    (\nabla_x f(\xbar_n,y_j)\zhat_n) \cdot\nabla_z u_0\\
    & = \int
    \mu_{\xbar_n}( dy_1) (\nabla_x f(\xbar_n,y_1)\zhat_n) \cdot\nabla_z u_0
    = (B_1(\xbar_n) \zhat_n ) \cdot\nabla_z u_0 \;,
  \end{aligned}
\end{equ}
 Likewise, using the independence of $Y^j_{x}$ for distinct $j$, the third term becomes 
\begin{equs}
\frac{1}{\lambda^2} \sum_{j,k=1}^\lambda & \int_0^\infty d\tau \BE \ftilde (\xbar_n, Y^j_{\xbar_n}(0))\otimes \ftilde (\xbar_n, Y^k_{\xbar_n}(\tau)) : \nabla_z^2 u_0 \\
& = \frac{1}{\lambda^2} \sum_{j=1}^\lambda  \int_0^\infty d\tau \BE \ftilde (\xbar_n, Y^j_{\xbar_n}(0))\otimes \ftilde (\xbar_n, Y^j_{\xbar_n}(\tau)) : \nabla_z^2 u_0 \\
 & = \frac{1}{\lambda} \int_0^\infty d\tau \BE \ftilde (\xbar_n, Y^1_{\xbar_n}(0))\otimes \ftilde (\xbar_n, Y^1_{\xbar_n}(\tau)) : \nabla_z^2 u_0 = \frac{1}{\lambda} \eta(\xbar_n) \eta(\xbar_n)^T : \nabla_z^2 u_0 \;.
\end{equs}
where the expectation is taken over realizations of $Y^j_{x}$ with $Y^j_{x}(0) \sim \mu_x$. Finally, since the $\nabla_{y_j} \BE_{y_k}$ term vanishes on the off-diagonal, the last term in \eqref{e:phmm_u0_big} reduces to
\begin{equs}
\frac{1}{\lambda}  \sum_{j,k=1}^\lambda & \int_0^\infty d\tau \int
\mu_{\xbar_n} (dy_j)\mu_{\xbar_n} (dy_k)  (\nabla_x b(\xbar_n,y_j)
\zhat_n)  \cdot \nabla_{y_j}
\BE_{y_k}\ftilde(\xbar_n,Y^k_{\xbar_n}(\tau))  \cdot \nabla_z u_0 \\ 
&= \frac{1}{\lambda}  \sum_{j=1}^\lambda  \int_0^\infty d\tau \int \mu_{\xbar_n} (dy_j)\mu_{\xbar_n} (dy_k)  (\nabla_x b(\xbar_n,y_j) \zhat_n)  \nabla_{y_j} \BE_{y_j}\ftilde(\xbar_n,Y^j_{\xbar_n}(\tau))  \nabla_z u_0 \\ 
& = \int_0^\infty d\tau \int \mu_{\xbar_n} (dy_1)\mu_{\xbar_n} (dy_k)
(\nabla_x b(x,y_1) \zhat_n)  \nabla_{y_1}
\BE_{y_1}\ftilde(\xbar_n,Y^1_{\xbar_n}(\tau))  \nabla_z u_0 \\
&= (B_2 (\xbar_n ) \zhat_n )  \cdot \nabla_{z} u_0 \;.
\end{equs}
It follows immediately that the reduced equation for the pair
$(\Xbar_n, \hat Z^\eps_n)$ is 
\begin{equs}
\frac{d\Xbar_{n}}{dt} &= F(\xbar_n) \\
d\Zhat_n &= B_0(\xbar_n) \Zhat_n dt + \lambda^{-1/2} \eta (\xbar_n) dV\;,
\end{equs}
with initial conditions $\Zhat_n(n\Delta t) = \zhat_{n}$ and $\Xbar_n(n\Delta t) = \xbar_n$.  Hence we see that $\zhat_{n+1}$ is described by 
\begin{equ}
\zhat_{n+1} = \zhat_n + B(\xbar_n) \zhat_n \Delta t +
\lambda^{-1/2}\eta(\xbar_n) \sqrt{\Delta t} \, \xi_n
\end{equ}
which is the Euler-Maruyama scheme for \eqref{e:phmm_clt_a}.

\subsection{Large fluctuations}
\label{s:phmm_ldp}

In this section we show that the LDP for PHMM is consistent with the
true LDP from Section \ref{s:fs_ldp}. In particular, let
\begin{equ}
u_{\lambda,\Delta t} (t,x) = \lim_{\eps \to 0} \eps \log \BE_x \exp \left(\eps^{-1} \vphi (x^\eps_{n}) \right)
\end{equ}
for $t \in I_{n,\Delta t}$, where $x^\eps_n$ is the PHMM approximation. We will argue that $u_{\lambda,\Delta t}(t,x) \to u(t,x)$ as $\Delta t \to 0$, where $u$ solves the correct Hamilton-Jacobi equation \eqref{e:fs_ldp_hj}.   
\par
The argument is a slight modification of that given in Section \ref{s:hmm_ldp}. Before proceeding, we recall the notation $S^{(\alpha)}_{\Delta t}$ for the semigroup associated with the Hamilton-Jacobi equation
\begin{equ}\label{e:phmm_ldp_hj}
\del_t u (t,x) = \CH (\alpha , \nabla u (t,x)) \;,
\end{equ}
where $\CH$ is the Hamiltonian defined by \eqref{e:fs_ldp_ham}. We also define the operator $S_{\Delta t} \vphi(x) = S^{(\alpha)}_{\Delta t} \vphi (x)|_{\alpha = x}$. 

As in Section \ref{s:hmm_ldp}, the claim follows from the asymptotic statement
\begin{equ}\label{e:phmm_ldp_dv}
\BE_x \exp \left( \eps^{-1} \vphi (x^\eps_n)  \right) \asymp \exp\left(\eps^{-1} (S_{\Delta t})^n \vphi(x) \right)\;, \quad \eps \to 0\;.
\end{equ}
Given \eqref{e:phmm_ldp_dv}, by an identical argument to that started in Equation \eqref{e:hmm_ldp_expansion}, it follows from \eqref{e:phmm_ldp_dv} that $u_{\lambda,\Delta t} $ is indeed a numerical approximation of the solution to \eqref{e:phmm_ldp_hj} and hence $u_{\lambda,\Delta t} \to u$ as $\Delta t \to 0$.   
\par
We will verify \eqref{e:phmm_ldp_dv} by induction, starting with the $n=1$ case. Since $(X^\eps_, \Ytilde^\eps_{0,1},\dots,\Ytilde^\eps_{0,\lambda})$ is a fast-slow system of the form \eqref{e:fs} with $\eps$ replaced by $\eps \lambda$, it follows from Section \ref{s:fs_ldp} (Varadhan's lemma) that 
\begin{equ}
\BE_x \exp \left((\eps \lambda)^{-1} \psi (X^\eps_1(\Delta t))\right) \asymp \exp \left( (\eps \lambda)^{-1} \Shat_{\Delta t}^{(\alpha)} \psi(x)|_{\alpha = x}  \right) \;,
\end{equ}
where $\Shat_{\Delta t}^{(\alpha)}$ is the semigroup associated with $\del_t v(t,x) = \CHhat (\alpha , \nabla v(t,x))$ and
\begin{equ}
\CHhat (\alpha,\theta) = \lim_{T\to \infty } T^{-1} \log \BE \exp \left(\theta \cdot \int_0^T d\tau \frac{1}{\lambda}\sum_{j=1}^\lambda f(\alpha, Y^j_{\alpha}(\tau)) \right)\;.
\end{equ}
Hence we have
\begin{equ}\label{e:phmm_ldp_semi}
  \begin{aligned}
    \BE_x \exp \left(\eps^{-1} \vphi (x^\eps_1)\right) &= \BE_x \exp
    \left((\eps \lambda)^{-1} \lambda\vphi (X^\eps_0(\Delta t))\right)\\
    &\asymp \BE_x \exp \left((\eps \lambda)^{-1}
      \Shat^{(\alpha)}_{\Delta t} (\lambda \vphi) (x)|_{\alpha = x}
    \right)
  \end{aligned}
\end{equ}
But since $Y^j_{\alpha}$ are iid for distinct $j$, the Hamiltonian $\CHhat$ reduces to
\begin{equs}
\lim_{T\to \infty } T^{-1} &\log \BE \exp \left(\theta \cdot \int_0^T d\tau \frac{1}{\lambda}\sum_{j=1}^\lambda f(\alpha, Y^j_{\alpha}(\tau)) \right)\\ 
&= \lambda \lim_{T\to \infty } T^{-1} \log \BE \exp \left(\frac{\theta}{\lambda} \cdot \int_0^T d\tau  f(\alpha, Y^1_{\alpha}(\tau)) \right) = \lambda \CH (\alpha , \frac{\theta}{\lambda} )\;.
\end{equs}
It follows that
\begin{equ}
\del_t \left(\lambda^{-1}\Shat_t^{(\alpha)} (\lambda \vphi)\right) = \lambda^{-1} \CHhat (\alpha , \nabla (\Shat_t^{(\alpha)} (\lambda \vphi))) = \CH (\alpha , \lambda^{-1}\nabla (\Shat_t^{(\alpha)} (\lambda \vphi)) )
\end{equ}
and hence $\lambda^{-1} \Shat^{(\alpha)}_{\Delta t} (\lambda \vphi) = S_{\Delta t}^{(\alpha)}\vphi$. Combining this with \eqref{e:phmm_ldp_semi} completes the claim for $n=1$. The proof of the inductive step for arbitrary $n\geq 1$ follows identically to Section \ref{s:hmm_ldp}.

\section{Numerical evidence}
\label{s:num}

In this section, we investigate the performance of the standard HMM
and PHMM methods for systems with well understood fluctuations and
metastability properties. These simple experiments confirm that HMM
amplifies fluctuations, which can drastically change the system's metastable
behavior, and that the PHMM succeeds in avoiding these problems. In
Section \ref{s:num_clt} we investigate simple CLT fluctuations for a
simple quadratic potential systems, in Section \ref{s:num_ldp} we look
at large deviation fluctuations for a quartic double-well
potential. Finally in Section \ref{s:num_nondiff} we look at
fluctuations for a non-diffusive double well potential, which has
large deviation properties that cannot be captured by a so-called
`small noise' diffusion.

\subsection{Small fluctuations} 
\label{s:num_clt}

We examine the small CLT-type fluctuations by looking the following
fast-slow system
\begin{equs}
\frac{dX}{dt} &= Y - X \\
dY &= \frac{\theta}{\eps} (\mu X - Y) dt + \frac{\sigma}{\sqrt{\eps}} dW\;.
\end{equs}
It is simple to check that the averaged system is given by
\begin{displaymath}
  \frac{d\Xbar}{dt} = (\mu - 1) \Xbar.
\end{displaymath}
Hence for $\mu < 1$ the
averaged system is a gradient flow in a quadratic potential centered
at the origin.

\begin{figure}[h]
\begin{center}
\includegraphics[width=0.5\textwidth]{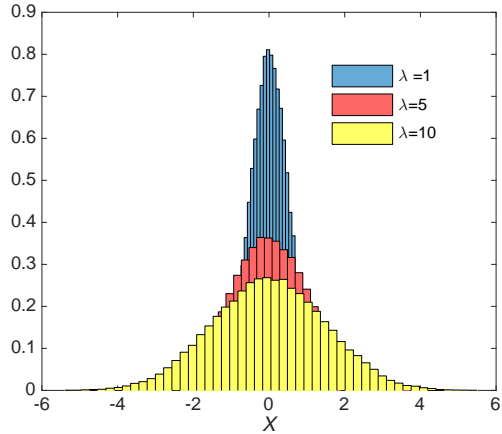}
\caption{Histogram of $X$ variables. Parameters used are
  $\eps = 10^{-2}$, $\delta t = 0.1$, $\theta= 1$, $\mu=0.5$,
  $\sigma= 5$, $T = 10^4$. }
\label{fig:clt_hist}
\end{center}
\end{figure}

\begin{figure}[h]
\begin{center}
\includegraphics[width=0.8\textwidth]{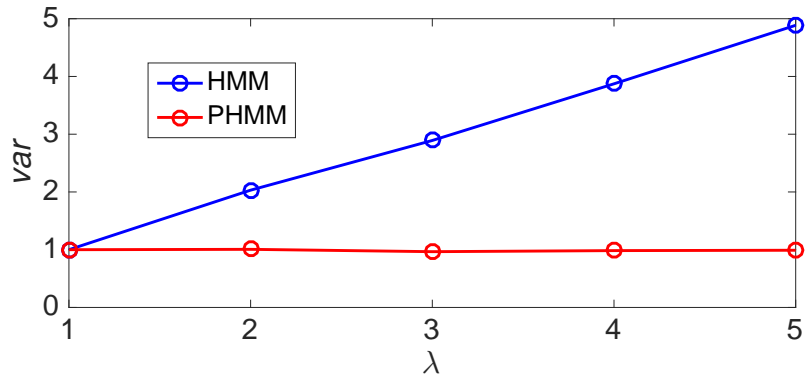}
\caption{Comparing the stationary variance of HMM and PHMM as a function of $\lambda$. Once again, we use parameters $\eps = 10^{-2}$, $\delta t = 0.1$, $\theta= 1$, $\mu=0.5$,  $\sigma= 5$, $T = 10^4$.}
\label{fig:clt_var}
\end{center}
\end{figure}

We will first illustrate that the HMM-type method described in Section
\ref{s:hmm} inflates the $O(\sqrt{\eps})$ fluctuations about the
average by a factor of $\sqrt{\lambda}$. In Figure \ref{fig:clt_hist}
we plot histograms of the slow variable $X$ for different speed-up
factors $\lambda$. It is clear that the spread of the invariant
distribution is increasing with $\lambda$. The profile remains
Gaussian but the variance is greatly inflated. In Figure
\ref{fig:clt_var} we plot the variance of the stationary time series
for $X$ as a function of $\lambda$. The blue line is computed using
HMM and the red line is computed using PHMM. As predicted by the
theory in Section \ref{s:hmm_clt}, in the case of HMM the variance is
increasing linearly with $\lambda$ and in the case of PHMM the
variance is approximately constant. Note that in this example, the CLT
captures the large deviations as well. 


\subsection{Large fluctuations}
\label{s:num_ldp}

To investigate the affect of parallelization on $O(1)$ deviations not
captured by the CLT, we will look at a fast-slow system which exhibits
metastability. Hence it is natural to take
\begin{equs}\label{e:num_fs2}
\frac{dX}{dt} &= Y - X^3 \\
dY &= \frac{\theta}{\eps} (\mu X - Y) dt + \frac{\sigma}{\sqrt{\eps}} dW\;.
\end{equs}
It is simple to check that the averaged system is
\begin{displaymath}
  \frac{d\Xbar}{dt} = \mu \Xbar - \Xbar^3.
\end{displaymath}
Hence for any $\mu > 0$ the
averaged system is a gradient flow in a symmetric double well
potential, with stable equilibria at $\pm \sqrt{\mu}$ and a saddle
point at the origin. The large fluctuations of the fast-slow system
can be investigated by looking at the first passage time for
transitions from a neighborhood of one stable equilibrium to the other.

\begin{figure}[h]
\begin{center}
\includegraphics[width=0.8\textwidth]{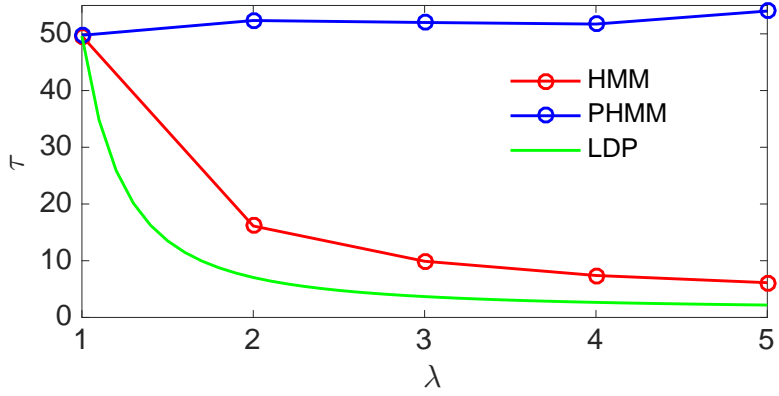}
\caption{The mean first passage time as a function of the speed-up
  factor $\lambda$, for HMM (red dotted) and PHMM (blue dotted). We
  include the LDP predicted curve for the mean first passage time of HMM, as
  discussed in Section~\ref{s:hmm_ldp}. $\eps = 10^{-3}$,
  $\delta t = 0.05$, $\theta= 1$, $\mu=1$, $\sigma= 15$,
  $T = 5\times 10^4$}
\label{fig:ldp_mfpt}
\end{center}
\end{figure}

In Figure \ref{fig:ldp_mfpt} we compare the mean first passage time
for HMM and PHMM as a function of $\lambda$. Even for $\lambda = 2$,
the distinction between the two methods is vast, with the mean first
passage time for HMM rapidly dropping off and for PHMM staying
approximately constant.

\begin{figure}[h]
\begin{center}
\includegraphics[width=\textwidth]{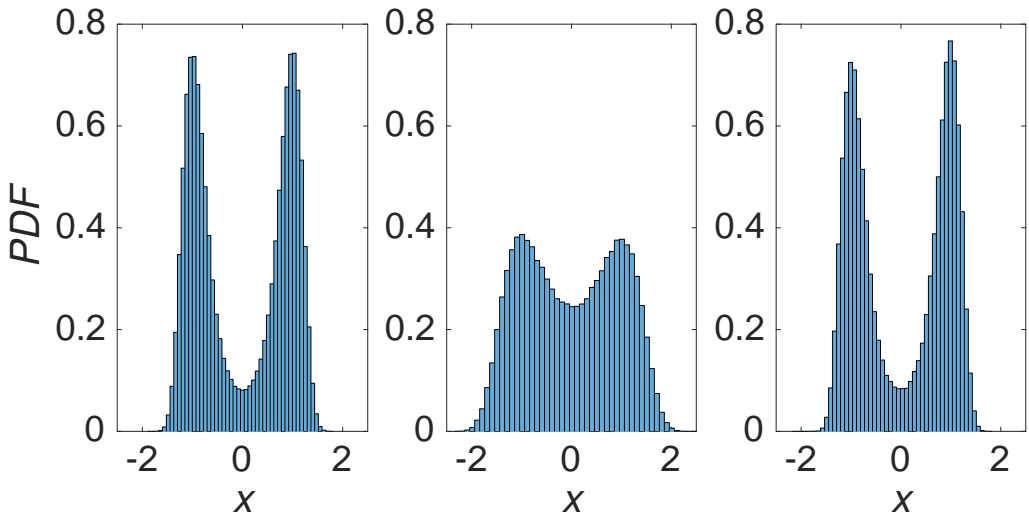}
\caption{Histogram of $X$ variables. $\eps = 10^{-3}$,
  $\delta t = 0.05$, $\theta= 1$, $\mu=1$, $\sigma= 15$,
  $T = 5\times 10^4$}
\label{fig:ldp_hist}
\end{center}
\end{figure}

In Figure \ref{fig:ldp_hist} we compare respectively the stationary
distributions of the true fast-slow system, HMM ($\lambda =5$) and
PHMM ($\lambda=5$). In the case of HMM, the energy barrier separating
the two metastable states is now overpopulated, which explains the
rapid fall in mean first passage time. In the case of PHMM, the
histogram is indistinguishable from the true stationary distribution
(with the exception of a slight asymmetry).

\begin{figure}
\centering
\includegraphics[width=0.8\textwidth]{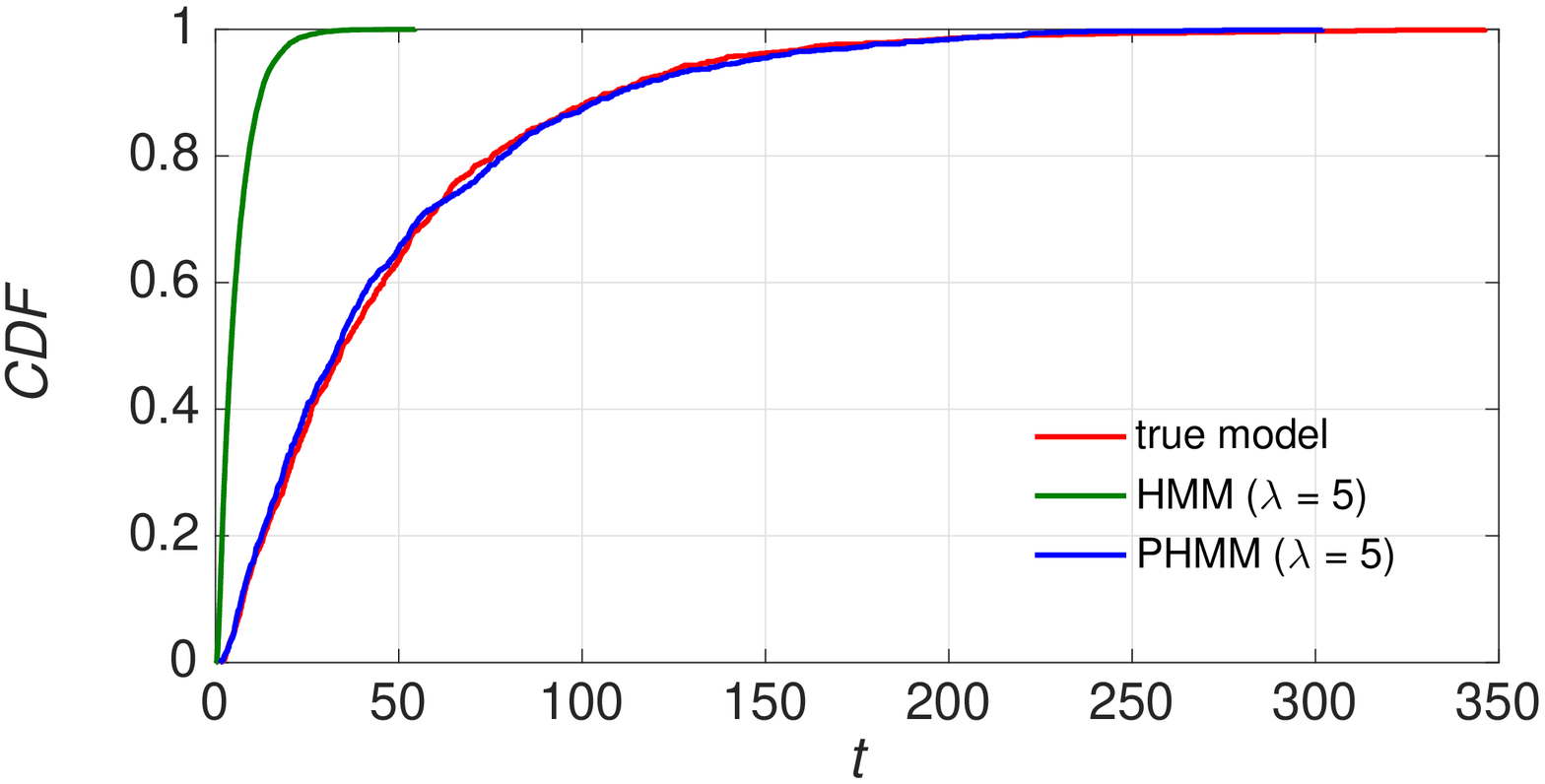}
\caption{Cumulative
  distribution functions for first passage times of the true model
  (red) for \eqref{e:num_fs2}, HMM with $\lambda = 5$ (green) and PHMM
  with $\lambda = 5$ (blue). The parameters used are $\eps = 10^{-3}$,
  $\delta t = 0.05$, $\theta= 1$, $\mu=1$, $\sigma= 15$,
  $T = 5\times 10^4$}\label{fig:ldp_fpt_cdf}
\end{figure}

In Figure \ref{fig:ldp_fpt_cdf} we plot the cumulative distributions
function (CDF) for the first passage time, comparing that of the true
fast-slow system, with HMM ($\lambda=5$) and PHMM ($\lambda=5$). We
see that the HMM first passage times are supported on a much faster
time scale than that of the true fast-slow system. In contrast, the
CDF of PHMM is almost indistinguishable from that of the true
fast-slow system. Hence PHMM is not just replicating the mean first
passage time, but also the entire distribution of first passage times.

\subsection{Asymmetric, non-diffusive fluctuations}
\label{s:num_nondiff}

\begin{figure}
  \centering
  \includegraphics[width=0.8\textwidth]{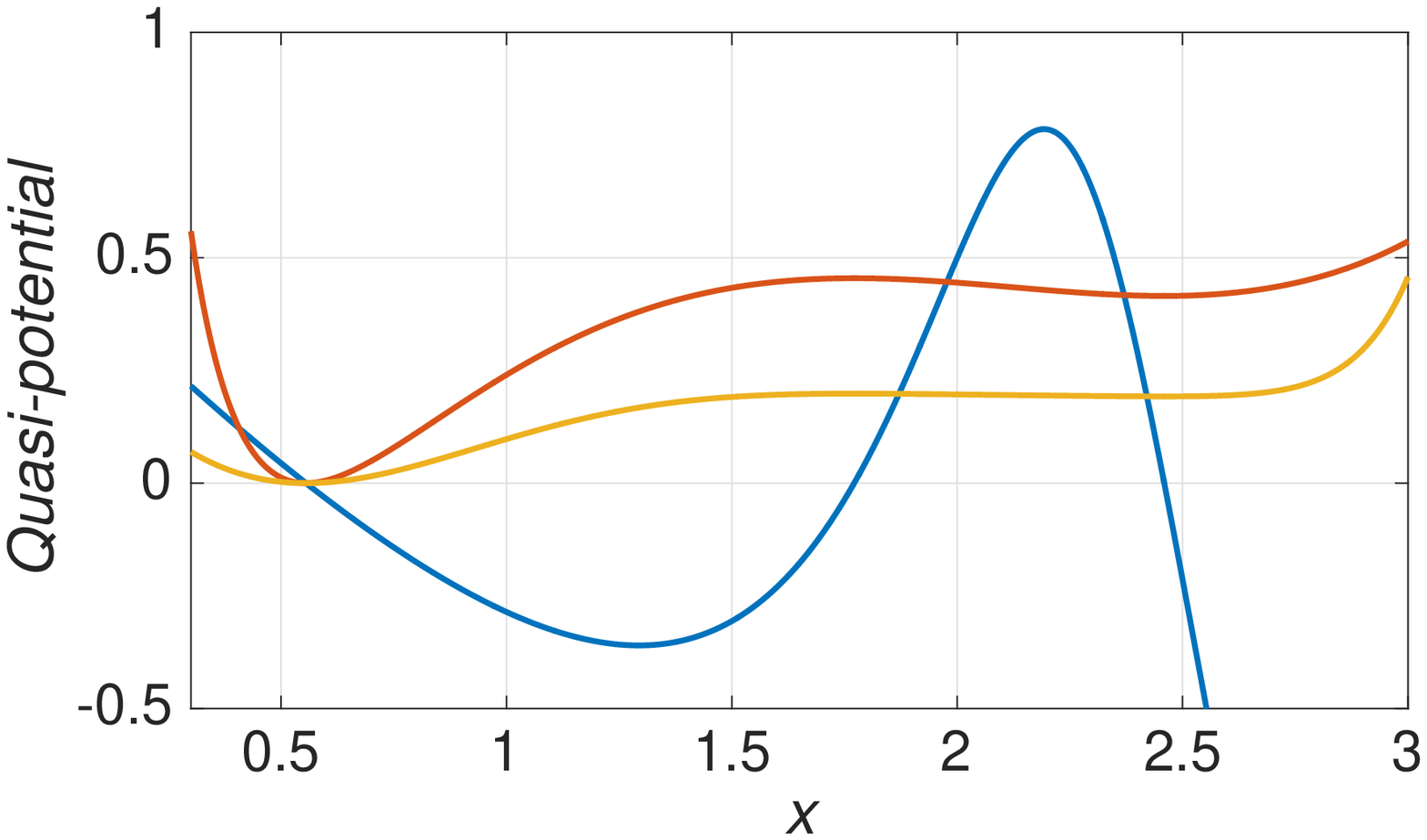}\caption{The quasi-potential
    $\CV(x)$ (red curve) and the one obtained from a quadratic approximation of
    the Hamiltonian (orange curve). Also shown in blue is the
    coefficient at the right-hand side of the reduced equation.\label{fig:num_quasi}}
\end{figure}

We now compare HMM and PHMM for a multiscale model that also displays
metastability, but in which the large fluctuations cannot be
characterized by a `small noise' Ito diffusion. In particular, the
Hamiltonian describing the LDP of the system is non-quadratic, as
opposed the the previous systems. The system has been used
\cite{bouchet15} to illustrate the ineffectiveness of diffusion-type
approximations for fast-slow systems. The fast-slow system is given by
\begin{equs}\label{e:num_fs3}
\frac{dX}{dt} &= Y^2 - \nu X\\
dY &=  - \frac{1}{\eps} \gamma(X) Y dt + \frac{\sigma}{\sqrt{\eps}} dW\;.
\end{equs}
where $\gamma(x) = x^4/10 - x^2 + 3$ .  The averaged equation for this
system reads
\begin{displaymath}
  \frac{d\bar X}{dt} = \frac{\sigma^2}{2\gamma(\bar X)} - \nu \bar X
\end{displaymath}
For $\nu = 1$ and $\sigma = \sqrt{3}$, this averaged equation
possesses two stable fixed points at $x\approx 0.555$ and
$x\approx=2.459$ and one unstable fixed point at $x\approx 2.459$. The
the rates of transition between these stable fixed points is captured
by the LDP. By an elementary calculation \cite{bouchet15}, the
Hamiltonian of this LDP is found to be non-quadratic and given by
\begin{displaymath}
  \CH (x,\theta) = - \nu x \theta +\frac12\left(\gamma(x) -
    \sqrt{\gamma^2(x)-2\sigma^2 \theta}\right)
\end{displaymath}
The quasi-potential associated with this Hamiltonian satisfies
$0= \CH(x,\CV')$, i.e.
\begin{displaymath}
  \CV'(x) = \frac{\nu x \gamma(x) - \tfrac12 \sigma^2}{\nu^2x^2}\;,
\end{displaymath}
and is displayed in Figure \ref{fig:num_quasi}. Whilst there is a
significant barrier corresponding to left-to-right transitions, there
is almost no barrier corresponding to right-to-left transitions.

\begin{figure}
\centering
\includegraphics[width=0.9\textwidth]{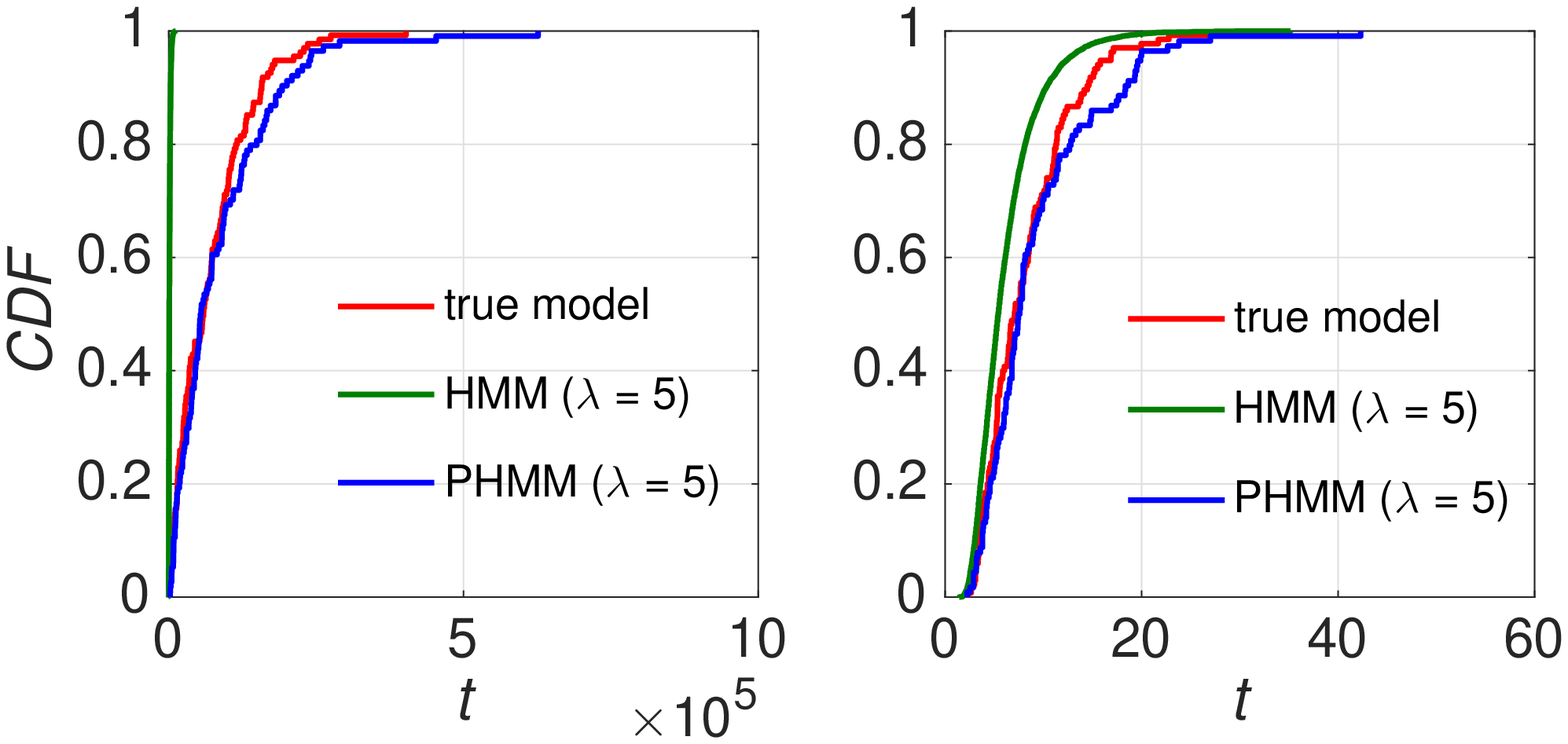}
\caption{Cumulative
  distribution functions for first passage times of the true model for
  \eqref{e:num_fs3} (red), HMM with $\lambda = 5$ (green) and PHMM
  with $\lambda = 5$ (blue). Left-to-right transitions on the left,
  right-to-left transition on the right.  The parameters used are
  $\eps = 0.05$, $\delta t = 0.5$, $\nu= 1$, $\sigma= \sqrt{3}$,
  $T = 1\times 10^7$}\label{fig:ldp_fpt_cdf2}
\end{figure}

In Figure \ref{fig:ldp_fpt_cdf2} we plot CDFs of the first passage
times: due to the asymmetry we plot separately the transitions from
the left-to-right and right-to-left.  For left-to-right transitions,
the HMM procedure drastically speeds up transitions because it 
enhances fluctuations: as is the case with the previous experiment,
the HMM transitions are supported on a timescale several orders of
magnitude faster than those of the true fast slow system. The PHMM
method does not experience this problem and the distribution of first
passage times agrees quite well with the true model.  For
right-to-left transitions, PHMM shows similarly good agreement with
the true fast-slow system, but in contrast HMM is not too far off
either. This can be accounted for by the `flatness' of the right
potential well, meaning that increasing the amplitude of fluctuations
will only decrease the escape time by a linear multiplicative
factor. We note that the noise appearing in the CDF plots is due to
the scarcity of transitions occurring in the model \eqref{e:num_fs3}.


\section{Discussion}
\label{s:discussion}

We have investigated HMM methods for fast-slow systems, in particular
their ability (or lack thereof) to capture fluctuations, both small
(CLT) and large (LDP). We found, both theoretically (Section
\ref{s:hmm_fluctuations}) and numerically (Section \ref{s:num}), that
the amplitude of fluctuations is enhanced by an HMM-type method. In
particular with an HMM speed up factor~$\lambda$, in the CLT the
variance of Gaussian fluctuations about the average is increased by a
factor~$\lambda$ as well. In the LDP, the quasi-potential is decreased
by a factor~$\lambda$, leading to the first passage times being
supported on a time scale $\lambda$ orders of magnitude smaller than
in the true fast slow system. This inability to correctly capture
fluctuations about the average suggests that HMM can be a poor
approximation of fast-slow systems, particularly when metastable
behavior is important. {As noted in Section
  \ref{s:hmm_ldp}, although the fluctuations of HMM are enhanced, the
  large deviation transition \emph{pathways} remain faithful to the
  true model. Thus we stress that HMM is a reliable method of finding
  transition pathways in metastable systems, but not for simulating
  their dynamics.}

We have introduced a simple modification of HMM, called parallel HMM
(PHMM), which avoids these fluctuation issues. In particular, the PHMM
method yields fluctuations that are consistent with the true fast slow
system for any speed up factor $\lambda$ (provided that we still have
$\eps \lambda \ll 1$), as was shown both theoretically (Section
\ref{s:phmm_fluctuations}) and numerically (Section \ref{s:num}). The
HMM method relies on computing one short burst of the fast variables,
and inferring the statistical behavior of the fast-variables by
extrapolating this short burst over a large time window. PHMM on the
other hand computes an ensemble of $\lambda$ short bursts, and infers
the statistics of the fast variables using the ensemble. Since the
ensemble members are independent, they can be computed in
parallel. Hence if one has $\lambda$ CPUs available, then the real
computational time required in PHMM is identical to that in HMM. 

Interestingly, one can draw connections between the parallel method
introduced here and the tau-leaping method used in stochastic chemical
kinetics \cite{gillespie00}. The tau-leaping method is an
approximation used to speed up simulation of stochastic fast-slow
systems of the type
\begin{equ}\label{e:tau_leap}
X^\eps(t) = X^\eps(0) + \sum_{k=1}^m \eps \CN_k 
\left(  \eps^{-1}\int_0^t a_k(X^\eps(s))ds  \right)\nu_k\;,
\end{equ}
where $\CN_k$ are independent unit rate Poisson processes, $\nu_k$ are
vectors in $\reals^d$ and $a_k : \reals^d \to \reals$. The system
\eqref{e:tau_leap} can be solved exactly by the stochastic simulation
algorithm (SSA), but when $\eps$ is small this can be extremely
expensive, due to the Poisson clocks being reset each time a jump
occurs. The tau-leaping procedure avoids this issue by chopping the
simulation window into sub-intervals of size $\tau$ and on each
subinterval fixing the Poisson clocks to their value at the left
endpoint. The speed-up is a result of the fact that one can simulate
the Poisson jumps in parallel, since their clocks are fixed over the
$\tau$ interval.
%
%
As a consequence of this analogy, one can check (using calculations
similar to those found above) that the tau-leaping method also
captures the fluctuations correctly, both at the level of the CLT and
that of the LDP. The former observation was made
in~\cite{anderson2011error}; to the best of our knowledge, the second
one is new.

As a final note, we stress that there are non-dissipative fast-slow
systems for which the PHMM will not be effective at capturing their
long time scale behavior, including metastability. These are system
for which the CLT and LDP hold on $O(1)$ timescale, but they either
cannot be extended to longer time-scale (in the case of the CLT) or
leads to trivial prediction on these time scales (in the case of the
LDP). To clarify this point, take for example the fast-slow
Langevin system
\begin{equs}
  \label{e:fs_hamiltonian}
  \dot{q}_1 &= p_1 \quad \dot{p}_1 = q_1 - q_1^3 + (q_2 - q_1) \;, \\
  \dot{q}_2 &= \eps^{-1} p_2 \quad \dot{p}_2 = \eps^{-1}(q_1 - q_2) 
- \eps^{-1} \gamma p_2 + \sqrt{2 \eps^{-1}\beta^{-1}\gamma} \eta \;.
\end{equs}
where $\gamma>0$ and $\beta>0$ are parameters.  For any value of
$\eps$, $\gamma$, this system is invariant with respect to the Gibbs
measure with Hamiltonian
\begin{equ}
  H(q_1,q_2,p_1,p_2) = \frac{1}{2}p_1^2 + \frac{1}{2}p_2^2 
  + \frac{1}{4}q_1^4 - \frac{1}{2}q_1^2 + \frac{1}{2}(q_1-q_2)^2 \;.
\end{equ}
As $\eps\to 0$, it is easy to check that the slow variables
$(q_1,q_2)$ converge to the averaged system
\begin{equ}
  \label{e:hamilt_avg}
  \Dot{\bar q}_1 = \bar p_1 \qquad \Dot{\bar p}_1 = - G'(\bar q_1)
\end{equ}
where the averaged vector field is the gradient of the free energy
\begin{equ}
  G(q_1) = \frac{1}{4}q_1^4 - \frac{1}{2}q_1^2 + \frac{1}{2}q_1^2 
= - \beta^{-1}\log \int \exp (\beta U(q_1,q_2))dq_2 \;,
\end{equ}
with
$U (q_1,q_2) = \frac{1}{4}q_1^4 - \frac{1}{2}q_1^2 +
\frac{1}{2}(q_1-q_2)^2$.
Likewise, if we introduce 
\begin{displaymath}
  \eta_1= \frac{q_1-\bar q_1}{\sqrt{\eps}}, \qquad 
  \zeta_1= \frac{p_1-\bar p_1}{\sqrt{\eps}},
\end{displaymath}
the CLT indicates that the evolution of these variables are captured
by 
\begin{equation}
  \label{eq:3}
  \dot \eta_1 = \zeta_1, \qquad d\zeta_1 = \sqrt{2\beta^{-1}\gamma}\, dB
\end{equation}
and we can also derive an LDP for \eqref{e:fs_hamiltonian} with
action
\begin{equation}
  \label{eq:4}
  \CS_{[0,T]} (q_1) = \frac{\beta}{4\gamma} \int_0^T |\ddot q_1-
    q_1+q_1^3|^2 dt
\end{equation}
However, neither~\eqref{eq:3} nor~\eqref{eq:4} capture the long
time behavior of the solution to~\eqref{e:fs_hamiltonian}. The problem
stems from the fact that the averaged equation in~\eqref{e:hamilt_avg}
is Hamiltonian, hence non-dissipative. As a result, fluctuations
accumulate as time goes on. Eventually, the CLT stops being valid, and
the LDP becomes trivial -- in particular, it is easy to see that the
quasi-potential associated with the action in~\eqref{eq:4} is
flat. For examples of this type, other techniques will have to be
employed to describe their long time behavior including, possibly,
their metastability (which, in the case of~\eqref{e:fs_hamiltonian} is
controlled by how small~$\beta^{-1}$ is, rather than~$\eps$). These
questions will be investigated elsewhere.

\bibliographystyle{./Martin}
\bibliography{./hmm_fluctuations}

\end{document}